
\input amstex
\input epsf
\documentstyle{amsppt}
\epsfverbosetrue
\NoBlackBoxes
\hyphenation{para-me-trised para-me-trise}

%
\define\fish{{1.5}} \define\whale{{1.12}} \define\crayfish{{1.13}}
\define\tick{{1.14}} \define\auk{{1.11}}
\define\shark{{2.1}} \define\parana{{2.2}} \define\ardvark{{3.6}}
\define\brontosaurus{{1.8}} \define\platypus{{1.9}} \define\frog{{4.3}}
\define\aligator{{4.5}} \define\panda{{5.1}} 
\define\lobster{{5.3}}  
\define\fly{{6.4}}
\define\antelope{{7.1}} \define\bandar{{8.3}}
\define\ibix{{12.4}} \define\emu{{12.5}}

%
\define\thesis{An}
\define\AB{AtBo}
\define\BHMM{BHMM}
\define\EGA {EGA}

\define\GH{GriHa}
\define\Hartshorne{Ha}
\define\Hi{Hi}
\define\La{La}
\define\murrayII{Mu}
\define\new{New}
\define\Po{Po}
\define\SaUhl{SaUhl}
\define\unitons{Uhl}
\define\Ward{Wd}
\define\Wr{Wr}
\define\Wooda{Wo}

%

\define\rock{\pi_{\R3}}\define\pebble{\pi_{T\P1}}

\define\p{\pmatrix}\define\pp{\endpmatrix}
\define\C{{\Bbb C}}

\define\tr{\operatorname{tr}}

\define\Un#1{{\operatorname{Harm}(\Bbb{S}^2,U(#1))}}
\define\Uns#1{{\operatorname{Harm}^*(\Bbb{S}^2,U(#1))}}

\define\U#1{{U}(#1)}
\redefine\u#1{\operatorname{u}(#1)}
\define\SU#1{{SU}(#1)}
\define\OO#1{{\Cal O}_{#1}}
\define\su#1{\operatorname{su}(#1)}
\define\gl#1{\operatorname{gl}(#1)}
\define\SLC#1{{SL}(#1,\C)}
\define\GLC#1{{GL}(#1,\C)}
\define\PGL#1{{PGL}(#1,\C)}
\define\slc#1{\operatorname{sl}(#1,\C)}
\define\glc#1{\operatorname{gl}(#1,\C)}
\redefine\O#1{\Cal O_{#1}} \define\I#1{\Cal I_{#1}} \define\E{\Cal E}
\define\HOM{\Cal Hom} 

\define\Aegt{{A^{g\gt}_\ebp}}\define\Aeg{{A^{g}_\ebp}}
  \define\CalAeg{{\Cal A^{g}_\ebp}}
   
\define\ue{{\Xi}}
\define\ul{{\Lambda}}

\define\Cs{{\Bbb C^\ast}}
\define\eb{{\bar\eta}}
\define\dbar{\bar\partial}
\define\ra{\to}
\define\ebp{{{\bar\eta}^\prime}}
\define\ep{{\eta^\prime}}
\define\eh{{\hat\eta}}\define\ehp{{\hat\eta{}^{\prime}}}
\define\eebu{{\pmatrix \eta\\ \bar{\eta}\\ u\endpmatrix}}
\define\F#1#2{F_{#1#2}}
\redefine\L#1#2{{L^{#1}_{#2}}}
\define\Gl#1{\operatorname{GL}(#1)}
\define\GL#1{\operatorname{GL}(#1)}
\define\Gr#1#2{{\text{Gr}_{#1,#2}}}
\define\pt{{\text{[pt]}}}

 \define\orth{{{}^{\text{th}}}}
\define\id{{\pmb 1}}
\define\norm#1{{\|#1\|}}
\define\Iff{{\it iff\ }}

\define\lp{{\lambda^\prime}} \define\lb{{\bar\lambda}}
\define\lbp{{\bar\lambda^\prime}} \define\lh{{\hat\lambda}}

\define\one{{\Bbb I}}

\redefine\P#1{{\Bbb P^{#1}}} \define\T{{\widetilde{T\Bbb P}{}^1}}
\define\R#1{\Bbb R^{#1}}
\define\th#1{\theta_{#1}} \define\ze#1{\zeta_{#1}}
\define\om#1{\omega_{#1}} \define\rh#1{\rho_{#1}}
\define\Rt#1{\tilde R_{#1}}
\define\st{{\tilde{s}}}

\define\toinfty{{\to\infty}}

\define\xh{{\hat{x}}}
\define\yh{{\hat{y}}}

\define\zhb{{\Hat{\bar{z}}}}\define\zbh{{\Hat{\bar{z}}}}
\define\zh{{\hat{z}}}
\define\zzbt{{\pmatrix z\\ \z\\ t \endpmatrix}}

\define\z{{\bar{z}}}
\define\deq{{\overset\text{def}\to=}}
\redefine\{{\left\lbrace} \redefine\}{\right\rbrace} \redefine\[{\left[}
\redefine\]{\right]} \redefine\({\left(} \redefine\){\right)}
\predefine\lef{\l}
\redefine\D#1{\nabla_{#1}}
\define\Du{\D u-i\Phi}

\redefine\dd#1{\frac\partial{\partial #1}}

\redefine\l{\lambda}
\define\gt{{\tilde g}} \define\gh{{\hat g}}
\define\Alg{A_\lb^g}\define\Algt{A_\lb^{g\gt}}\define\Algh{A_\lb^{g\gt\gh}}
\define\dB#1{{\dsize\frac{\partial B}{\partial#1}{\ssize\zzbt}}}
\define\eval{\operatorname{eval}}

\define\figa{{Fig.\,1}} \define\figg{{Fig.\,2}} \define\figj{{Fig.\,3}}
\define\figh{{Fig.\,4}}

\define\apict{{\midinsert
\hfill\epsfxsize=.7\hsize\epsffile{a.epsf}\hfill
\botcaption{\figa}
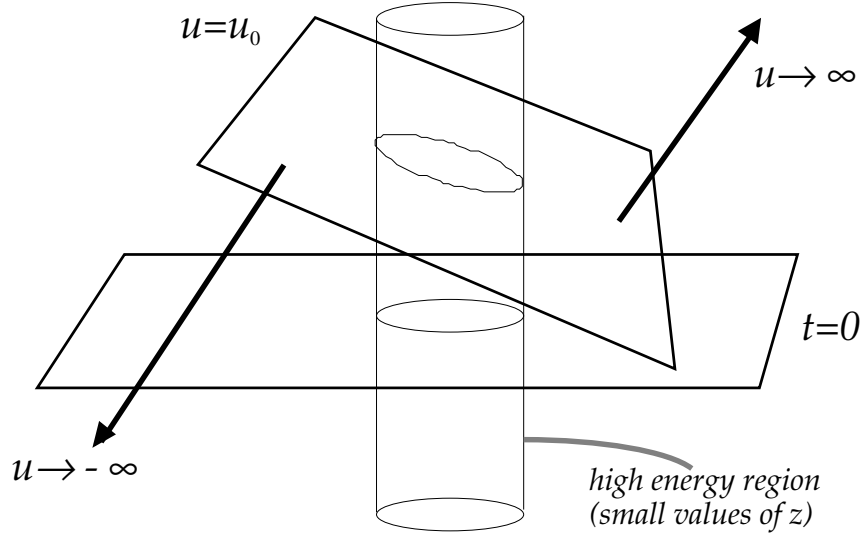
A picture of $\R3$ showing the high energy cylinder and a plane
$\{u=u_0\}$ for
some direction $\l\in\C^*, |\l|\ne1$.
{}From the picture, we would expect solutions to $\D\ebp$ to extend to
$\eta=\infty$, and solutions of $\D u-i\Phi$ to extend to $u=\infty$,
since in
these directions the connection coefficients decrease rapidly.
\endcaption
\endinsert  }}

\define\gpict{{\midinsert
\hfill\epsfxsize=.8\hsize\epsffile{g.epsf}\hfill
\botcaption{\figg}
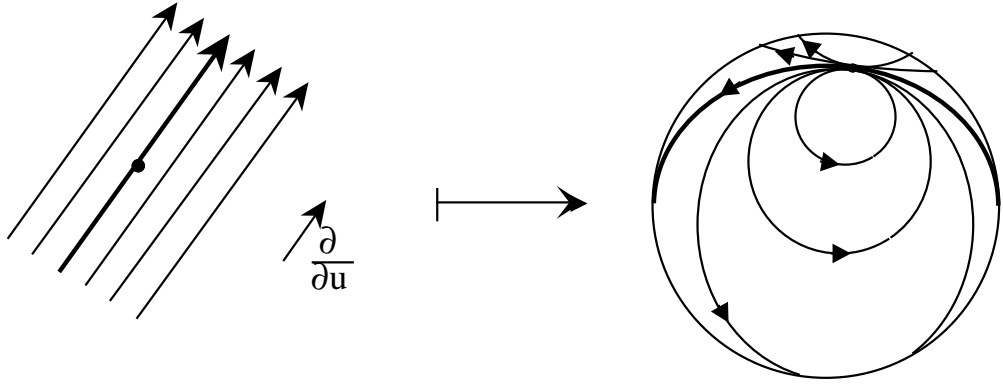
For a fixed value of $\l\in\Cs$, lines parallel to $\dd u$ in $\R3$ are
projected onto lines on $\R2$ which are completed as circles on ${\Bbb
S}^2$
tangent at the point at infinity.
\endcaption
\endinsert  }}

\define\hpict{{\midinsert
\hfill\epsfxsize=.8\hsize\epsffile{h.epsf}\hfill
\botcaption{\figh}
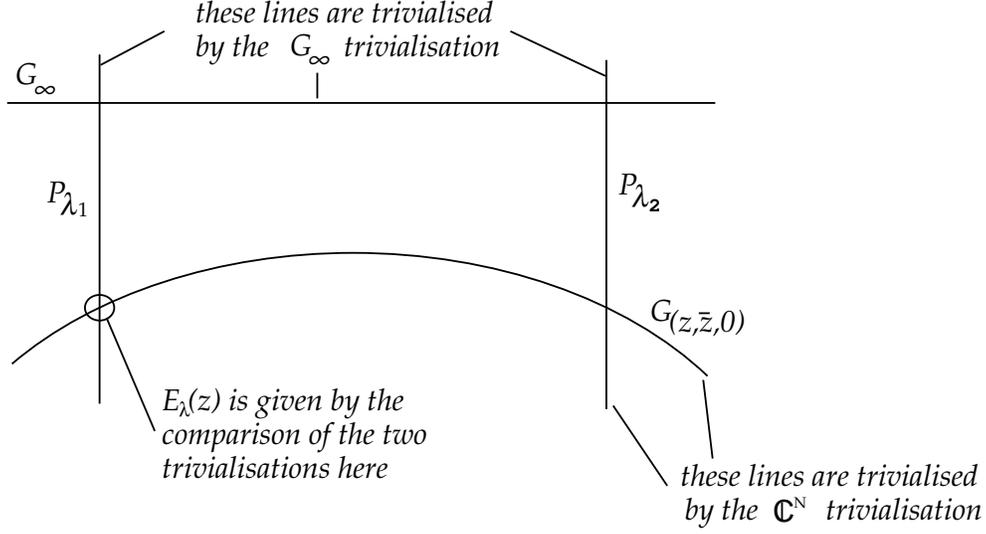
Compare the trivialisations at their point of intersection.
\endcaption
\endinsert  }}

\define\jpict{{\midinsert
\hfill\epsfxsize=.8\hsize\epsffile{j.epsf}\hfill
\botcaption{\figj}
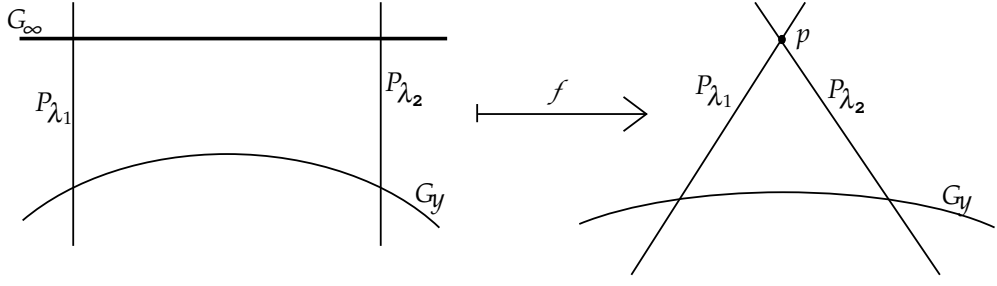
The embedding of $\T\to Q\subset\P3$ maps the section at infinity to a
singular point.  Hyperplane sections of $Q$ pull back to sections,
${G}_y$,
of $T\P1$, or to unions $P_{\l_1}\cup {G}_\infty\cup P_{\l_2}$ if they
contain the singular point.
\endcaption
\endinsert  }}

\topmatter
\title Uniton Bundles \ \ dg-ga/9508011\endtitle
\author Christopher Kumar Anand \endauthor
\affil McGill University\endaffil
\address \noindent Department of Mathematics and Statistics,
 McGill University,
  Burnside Hall,\linebreak 805 ouest, rue Sherbrooke,
 Montr\'eal\ \ QC\ \  H3A 2K6\endaddress
\curraddr Mathematics Research Centre,
  University of Warwick,
  Coventry CV4 7AL\endcurraddr
\email anand\@maths.warwick.ac.uk \endemail
\date 18 July 1994\quad  Revised:\ 4 August 1995\enddate
\thanks Research supported by NSERC and FCAR scholarships. \endthanks
\abstract
We show that a twistor construction of Hitchin and Ward can be adapted to
study unitons (harmonic spheres in a unitary group).  Specifically, we
show that unitons are equivalent to holomorphic bundles with extra
structure over a rational ruled surface with energy given by Chern class.
This equivalence allows us to
confirm the conjecture of Wood that unitons are rational. We also
construct an example of a two uniton in U(3) using this construction.
\endabstract
\toc \widestnumber\head{19.}
\head 1. Introduction\page 1\endhead
\head 2. Bogomolny Equations\page4\endhead
\head 3. Twistors\page5\endhead
\head 4. Adapted Coordinates\page{7}\endhead
\head 5. Compactness\page{9}\endhead
\head 6. Extra Structure\page{15}\endhead
\head 7. $\T$ as a conic\page{16}\endhead{}
\head 8. Compact Twistor Fibration\page{18}\endhead{}
\head 9. The Connection and Higgs' Field \page{19} \endhead
\head 10. Extra structure\page{20}\endhead
\head 11. Chern class equals energy\page{22} \endhead
\head 12. Ward's Construction\page{24}\endhead
\head 13. Example of a U(3) uniton\page{26} \endhead
\head 14. Wood's Conjecture\page{28}\endhead
\endtoc
\endtopmatter

\document
\pageheight{9 true in} 
\pagewidth{6.5 true in}

\head 1. Introduction\endhead
Unitons are harmonic maps  $S:{\Bbb S}^2\to U(N)$, that is,
maps
satisfying
$$\dd x(S^{-1}\dd x S)+\dd y(S^{-1}\dd y S)=0.\tag1.1$$
More generally, harmonic maps between Riemannian manifolds $M$ and $N$
are critical values of an energy functional
$$\text{energy}(\phi:M\to N)=\frac12\int_M|d\phi|^2.$$
In the case of maps into a matrix group, with
the standard (left-invariant) metric, the energy takes the form
$$\text{energy}(S)=\frac12\int_{\R2}
\left(|S^{-1}\dd xS|^2+|S^{-1}\dd y S|^2\right) dx\wedge
dy.\tag1.2$$
The uniton equations are the corresponding Euler-Lagrange equations.

{}From \cite{\SaUhl, Theorem 3.6}, we know that
harmonic maps from $\R2\to \U N$ extend to ${\Bbb S}^2$ \Iff they have
finite
energy, and that such maps are always smooth.
In the following, we will use this fact and work in terms of coordinates
$x$ and $y$ on $\R2$.

Unitons are determined by $A=A_xdx+A_ydy$ ($=A_zdz+A_\z d\z$ in
characteristic coordinates)
$$A_x\deq\frac12S^{-1}\dd x S,\quad A_y\deq\frac12S^{-1}\dd y S,\tag1.3$$
and a choice of initial condition, $S(\infty)\in\U N$, as we can see by
thinking
of $d+2A$ as a flat connection, and $S$ as a gauge transformation.
(The choice of initial condition or basing condition gives a covariant
constant basis of the bundle.)  We thus have a decomposition of the
 $U(N)$-unitons $\Un N$:
$$\Un N=\U N\times\Uns N,\tag1.4$$
where
$$\Uns N\deq\{S\in\Un N:S(\infty)=1\}$$
will be called the based unitons.  Of course, the energy doesn't depend
on the basing condition, and we can write it in terms of $A$ as
$$\text{energy}=-4i\int \operatorname{tr} A_zA_\z\,dz\,d\bar{z}.$$

Since $\Uns N=\{A=dS:S\in\Un N\}$, it is useful to
have equations for $A$ as well.
Two $\u N$-valued maps $A_x,A_y$ come from a map $S:\R2\to\U N$ in this
way \Iff $d+2A$ has zero curvature ($S$ is the flat gauge) \Iff
$$\align0&=d(2A)+[2A,2A]=2\{\dd x A_y-\dd y A_x+2[A_x,A_y]\}dx\wedge
dy.\tag\fish\ a\\
\intertext{They come from a {\it harmonic\/} map if in addition}
0&=d^*A=\dd x A_x+\dd y A_y.\tag\fish\ b\endalign$$
    The map $S:\R2\to\U N$ extends to a smooth map ${\Bbb S}^2\to\U N$ \Iff
$$A_\zh\deq-z^2A_z,\text{ and }A_\zhb\deq-\z^2A_\z\tag1.6$$
are smooth at
$z=\infty$, where we make use of complex coordinates $z=x+iy,\zh=1/z$.
In terms of complex coordinates, the uniton equations are (any two of)
$$\aligned\dd \z A_z-\dd z A_\z+2[A_\z,A_z]&=0,\\ \dd \z A_z+[A_\z,A_z]
&=0,\\
\dd z A_\z+[A_z,A_\z] &=0.\endaligned\tag1.7$$

In this paper we prove
\proclaim{Theorem A}  The space of based unitons, $\Uns N$, is
isomorphic to the
space of  $U(N)$-uniton bundles, with energy corresponding to the second
Chern class.\endproclaim
\noindent
where the  uniton bundles are bundles on $\T\deq \P{}(\O{}\oplus\O{}(2))$,
the fibrewise compactification of the tangent bundle $T\Bbb P^{1}$ of the
{\it complex}
projective line.

Let $(\l,\eta)$  and
$(\lh=1/\l,\eh=\eta/\l^2)$ be coordinates on $T\P1\cong \OO{\P1} (2)$,
where $\l$ is the usual coordinate on $\P1$ and $\eta$
is the coordinate associated to $d/d\l$.
Meromorphic sections ($s$) of $T\P1$ give all the holomorphic sections of
$\T$ ([s,1] in projective coordinates), except the section at infinity
($[1,0]$).
We will use the following notation for curves on $\T$:
$$\aligned P_\l&=\pi^{-1}(\l\in\P1)=\text{a pfibre (silent p)}\\
{G}_0&=\{(\l,[0,1])\}=\text{(graph of) zero section of $T\P1$}\\
{G}_\infty&=\{(\l,[1,0])\}=\text{infinity section of $\T$}\\
{G}_{\eta=s}&=\{(\l,[s(\l),1])\}.\endaligned\tag\brontosaurus$$
If $y=(a,b,c)\in\C^3$, we will also write ${G}_y$ for
${G}_{\eta=\frac12(a-2b\l-c\l^2)}$.

We define two  \lq  real structures\rq, i.e.\ antiholomorphic involutions,
 on $T\Bbb P^{1}$ by
$$\tau^*(\l,\eta)=(-1/\lb,-\lb^{-2}\eb),\quad
\sigma^*(\l,\eta)=(1/\lb,-\lb^{-2}\eb).\tag\platypus$$
On $\C^3\cong H^{0}(\Bbb P^{1},\Cal O(2))$, the space of finite sections,
 they act by
$$\tau^*(a,b,c)=(\bar c,\bar b,\bar a),\quad \sigma^*(a,b,c)=(\bar c,-\bar
b,\bar a).$$
Since these real structures are related by time translation
$$\align\delta_{t}:\ & (\lambda,\eta)\mapsto (\lambda,\eta -2t\lambda)\\ &
(a,b,c)\mapsto (a,b+t,c)\endalign$$
(real time for one is imaginary time for the other) we need only one real
structure $(\sigma)$.  We use both to simplify the exposition and to
maintain contact with previous results for monopole bundles which we will
require.

\proclaim{Definition of Uniton Bundles} A rank $N$, or $\U N$,
{\rm uniton bundle},
$\Cal{V}$,
is a holomorphic rank $N$ bundle on $\T$ which is a) trivial when restricted
to the following curves
in $\T$
\roster
\item the section at infinity
\item nonpolar fibres (i.e.\ fibres above $\l\in\Cs\subset\P1$)
\item real sections of $T\P1$
\endroster
b) is equipped with bundle lifts
$$\CD \Cal{V}@>\tilde\delta_t>>\Cal{V}\\@VVV
@VVV\\\T@>\delta_t>>\T\endCD\quad\text{and}\quad
\CD \Cal{V}@>\tilde\sigma>>\Cal{V}^*\\@VVV @VVV\\\T@>\sigma>>\T\endCD$$
\roster
\item $\tilde\delta_t$, a one-parameter family of holomorphic
transformations fixing $\Cal{V}$ above the section at infinity,
lifting $\delta_t$  and
\item $\tilde\sigma$ a norm-preserving, antiholomorphic lift of
 $\sigma$ such that the induced
hermitian metric on $\Cal{V}$ restricted to a fixed point of $\sigma$
is positive definite,
equivalently, the induced lift to  the
principal bundle of frames  acts on fibres of fixed points of $\sigma$ by
$X\mapsto X^{*-1}$.\endroster
and c) has a framing, $\phi\in H^0(P_{-1},Fr(\Cal{V}))$, of the bundle
$\Cal{V}$ restricted
to the fibre $P_{-1}=\{\lambda=-1\}\subset \T$ such that
$\tilde\sigma(\phi)=\phi$.
\endproclaim

Theorem A
 extends a result of Ward \cite\Ward\ that harmonic maps $\R2\to\SU N$
correspond to bundles on $T\P1$ and that when $N=2$, they extend to
bundles on the compactified space $\T$. It also implies that the energy
spectrum is discrete, as shown in \cite{Va1}.

The space $\T$ is a rational ruled surface, and we will see from the
existence of the time translation $\tilde{\delta}_{t}$
that the bundle $\Cal{V}$ restricted to
fibres, $P_\l$ of $\T\to\P1$, is trivial except for $P_0$ and $P_\infty$
which are \lq  jumping lines\rq\ (\cite{{Hu1}}).
By proving Theorem A we can study the topology of $\Uns N$ via
the spaces of allowed jumps at $P_0$ and $P_\infty$ and of monads.
These are the techniques which
were
successful in proving the Atiyah-Jones Conjecture for Yang-Mills
instantons (see \cite\BHMM).  Of course, this line of thinking leads to
many new questions.  Is the uniton number tied to the holomorphic
structure of the bundle near $P_0$ and $P_\infty$? Is there a
closed-form expression for
the uniton in terms of monad data?

These questions will be the subject of a future paper.  In this paper
we restrict ourselves to one example and one straightforward
application: a
reaffirmation of Wood's conjecture that unitons are composed of rational
functions of $x,y\in\R2$:
\proclaim{Corollary B}
If $S:{\Bbb S}^2\to \U N$ is a uniton, then the composition with
$\U N\hookrightarrow\Gl N$
is rational, i.e.\ the functions in $x$ and $y$ which make up the matrix
$S\in \U N$ are rational.
\endproclaim
As pointed out to the author by Martin Guest and the referee,
Wood's Conjecture was
affirmed by Valli \cite{Va2} by showing that the lift of the harmonic map to
the loop group (i.e.\ the extended solution) is real algebraic.  Valli's use
of interpolation theory is technically quite different from our treatment,
and his proof applies to the more general case of pluriharmonic maps
from simply-connected compact complex projective manifolds.
For a discussion of the broader historical development and current
motivation
of this problem
see \cite\thesis.

To prove Theorem A we will extend a twistor construction of Hitchin and
Ward (reviewed in \S3).  This is made possible by embedding the unitons
in the solutions of the Bogomolny equations (\S2).  Both geometric
argument and Sobolev methods (\S5) are involved in extending the
constructed
bundle to the compactified space.  Algebro-geometric methods are
required to show that the inverse construction goes through (\S\S7-10).
We complete the proof of Theorem A in \S11 by computing the Chern-Weil
integral.
In \S12 we relate our construction to  Ward's construction, by means of
which
we construct a $U(3)$ uniton in \S13 and  affirm Wood's conjecture
in \S14.

\subhead 1.10 Extended Solutions\endsubhead
We will also make use of Uhlenbeck's extended solutions $E_\l$ (actually
first employed in \cite\Po), which encode the unitons as follows
\proclaim{Theorem \auk\ \cite{\unitons, 2.1}}  Let $\Omega\subset {\Bbb
S}^2$ be a simply-connected
neighbourhood and
$A:\Omega\ra T^*(\Omega)\otimes U(N)$.
Then $2A=S^{-1}d S$, with $S$ harmonic \Iff the curvature of the
connection $$\Cal D_\l=(\dd \z+(1+\l)A_\z,\dd z+(1+\l^{-1})A_z)\tag\whale$$
vanishes for all $\l\in\C^*$.\endproclaim
\proclaim{Theorem \crayfish\ \cite{\unitons, 2.2}} If $S$ is harmonic and
$S(\infty)=\one$, then there exists a unique covariant constant
frame $E_\l:\P1\to
\U N$ for the connection $\Cal D_\l$ for $\l\in\C^*$ with
\roster \item $E_{-1}=\one$, \item $E_1=S$,
\item
$E_\l(\infty)=\one$.\endroster
  Moreover, $E_\l$ is analytic and holomorphic in
$\l\in\Cs$.\endproclaim
\proclaim{Theorem \tick\ \cite{\unitons, 2.3}}
Suppose $E:\Cs\times\Omega\to G$ is
analytic and holomorphic in the first variable, $E_{-1}\equiv\one$, and
the expressions
$$\frac{E_\l^{-1}\dbar E_\l}{1+\l},\quad\frac{E_\l^{-1}\partial
E_\l}{1+\l^{-1}}$$
are constant in $\l$ then $S=E_1$ is harmonic.\endproclaim
The extended solution is the leitmotiv of this paper.  It encodes the
analytic properties of the uniton and without it we would be forced to
use a lot more analysis in proving the compactness result of Theorem A.
Moreover, the energy calculation, the construction of a U(3) uniton, and
proof of Wood's conjecture use it
essentially.

\head 2. Bogomolny Equations\endhead

The next important rewriting of the uniton equations was
Ward's embedding of the harmonic map equations into the
Bogomolny equations over $\Bbb R^{2,1}$ as time-independent
 solutions \cite\Ward.
Since the geometry is simpler, we prefer to work with the Euclidean
equations:
 $$\D{}\Phi=*F,$$
where $\D{}=d+A$ is a connection, $\Phi$ is a Higgs' field, the curvature
$F=\D{}\circ\D{}=dA+A\wedge A$ and the Hodge-star is given by $*d y\wedge
d t=d x$, $*d t \wedge d
x=d y$,
and $*d x\wedge d y= d t$.
Assuming time independence of $\D{}$ and $\Phi$, the equations are
$$\aligned [A_t,\Phi]&=\dd x A_y-\dd y A_x+[A_x,A_y],\\
\dd x \Phi+[A_x,\Phi]&=\dd y A_t+[A_y,A_t],\\
\dd y \Phi+[A_y,\Phi]&=-\dd x A_t+[A_t,A_x].\endaligned\tag\shark$$
Write the system \thetag\fish\ as $\D x A_x=-\D y A_y$ and $\D x A_y=\D y
A_x$, and
\thetag\shark\ as $\D x A_y-\D y A_x-[A_x,A_y]= [A_t,\Phi]$, $\D
x\Phi=\D y A_t$, and $\D y \Phi=-\D x A_t$.
If $$A_t=-iA_y,
\ \ \ \Phi=iA_x,\tag\parana$$ the two systems are equivalent.
It is important to note that $A_t$ and $\Phi$ are imaginary (i.e.\ in $i\u
N$),
as this will determine the real structure we will use on $T\P1$.

We can use the freedom to change gauge to put any
$t$-independent solution
into the form \thetag{\parana}; the new gauge, $g$, is
given by solving $\dd y g\, g^{-1}=iA_t-A_y,\ \dd x g\,g^{-1}=-i\Phi-A_x$,
which we can do because the appropriate curvature component
$$\align\[\dd x+A_x+i\Phi,\dd y+A_y-iA_t\]&=\(\dd x A_y-\dd y
A_x+[A_x,A_y]-[A_t,\Phi]\)\\
&\ \ \ \ +i\(-(\dd x A_t+[A_x,A_t])-(\dd y\Phi+[A_y,\Phi])\)\\
&=\(F_{xy}-\D t\Phi\) +i\(F_{tx}-\D y \Phi\)\endalign$$
vanishes for solutions.  Thus a uniton is equivalent to a time
independent solution of the Bogomolny equations.
Of course, we still have to solve $S^{-1}dS=2A$ (i.e.\ integrate) to get a
uniton.

For future reference, we extract from the previous discussion the
following
\proclaim{Theorem 2.3} The space of based unitons, $\Uns N,$ is
isomorphic
to the space of $t$-independent solutions $\{(\D{},\Phi)\}$ to the
Bogomolny equations with finite
energy $\int_{\R2}|A_x|^2+|A_y|^2<\infty$ (equivalently, such that
$\lim_{x+iy\to\infty}A_z/z^2$ exists), which are real in the sense that
$A_x,A_y\in\u N$ in a gauge such that $A_t=-iAy$ and $\Phi=iA_x$.
\endproclaim

\head 3. Twistors\endhead

Oriented lines in $\R3$ are given by a direction and a displacement from
the origin perpendicular to the line's direction.  Collectively, they make
up the space $T {\Bbb S}^2\cong T\P1=\OO{\P1}(2)$ (the exact
correspondence
depending
on an isomorphism of ${\Bbb S}^2$ and $\P1$ to be specified using
stereographic projection).  Reversing the direction of geodesics
corresponds to the antiholomorphic involution, $\tau$, of $T\P1$ which is
the
negative of the map on $T\P1$ induced by the antipodal map,
$\tau^*\l=-1/\lb$, on $\P1$.  The holomorphic bundle $T\P1$ has sections
$\eta=a+b\l+c\l^2$, where $\l$ and $\eta$ are base and fibre coordinates
over $\P1\setminus\{\infty\}$.  So its section space is $\C^3$.

Let $\underline{\R{}}$ be the trivial real line bundle.
The geometry can be represented by a {\it twistor fibration}:
$$\matrix
 & & {\Bbb S}^2\times\R3\cong T\P1\oplus\underline{\R{}}\\\phantom{H}\\
 \rock & \swarrow &  &\searrow&\pebble\\\phantom{H}\\
\R3& & & & T\P1.\endmatrix\tag3.1$$
The explicit twistor correspondence, associating points of $T\P1$
to null planes
in $\C^3$ and points of $\C^3$ to sections of $T\P1$, is
$$\aligned (\l,\eta)\in T\P1&\mapsto
\{(a,b,c)\in\C^3:\eta=\frac12(a-2b\l-c\l^2)\},\\
 (a,b,c)\in\C^3&\mapsto \{(\eta,\l)\in
T\P1:\eta=\frac12(a-2b\l-c\l^2)\},\endaligned\tag3.2$$
where $\R3\hookrightarrow\C^3$ as
$$(x,y,t)\mapsto(x+iy,2t,x-iy).$$

The point of this construction is that there
is a {\it twistor correspondence} between solutions $(\D{},\Phi)$ to the
Bogomolny equations on $\R3$ and holomorphic bundles on $T\P1$ which are
trivial on {\it real sections} of $T\P1$.  This twistor construction is due
to Hitchin and Ward.  See \cite{{Hi}} for details.

\subhead 3.3 The Bundle\endsubhead
Let $\Cal{W}=\C^N\times\R3$ be the trivial bundle over $\R 3$.  Define the
bundle $\Cal{V}\to T \Bbb P^1$, over $\ell\in T \Bbb P^1$ by
$$\Cal{V}_{\ell}=\{s\in H^0(\ell,\rock^*\Cal{W})\,:(\D
u-i\Phi)s=0\}.\tag3.4$$
where the line $\ell\subset\R3$ is parametrised by
arclength, $u$, and $\nabla_{u}$ represents the covariant derivative
corresponding to $\dd u$.

The bundle $ \Cal{V}$ comes with a natural $\dbar$-operator, i.e.\ an
operator
$$\dbar:\Gamma(T\P1, \Cal{V}\otimes T^{(p,q)}T\P1))\ra \Gamma(T\P1,
\Cal{V}\otimes
T^{(p,q+1)}T\P1),$$
which satisfies a Leibnitz rule and $\dbar^2=0$.
This defines a complex bundle structure (since its flat sections will give
local holomorphic framings).

The $\dbar$-operator is defined as follows.

\roster
\item Pull back the connection $(\Cal W, \nabla)$ on $\Bbb R^{3}$
to $(\pi^{*}\Cal W,\pi^{*}\nabla)$ by the projection $\pi :\Bbb
R^{3}\times\Bbb P^{1}\to\Bbb R^{3}$.
\item The embedding $i:\P1\hookrightarrow\R3$ induces a splitting of
$\R3\times\P1\to\P1$ into $$T\P1\oplus<\text{normal bundle}>.$$\endroster
On the first component, $T\P1$,
we put the restriction of $\pi_{\R3}^*\D{}$ on $\Bbb R^3\times\Bbb P^1$.
  On the normal bundle, we
put $\pi_{\R3}^*(\D u-i\Phi)$.

Since $T\P1$ has a complex
structure, we can split $T(T\P1)$ into holomorphic and antiholomorphic
parts.  The restriction of the connection to the antiholomorphic part,
$\overline{\nabla}$, acts on sections of $\Cal{W}|_\ell$ and since
$\{\overline{\nabla},\D
u-i\Phi\}$ are involutive (as a consequence of the Bogomolny equations)
it induces a $\dbar$-operator on $\Cal{V}$.

\proclaim{Theorem 3.5 \cite{\murrayII, 18}}
If $(\D{},\Phi)$ is a unitary solution of  the Bogomolny equations
$\D{}\Phi=*F$ on $\R3$, then $\Cal{V}$ is in a natural way a holomorphic
bundle on the space of geodesics $T\P1$ such that
\roster
\item $\Cal{V}$ is trivial on every real ($\tau$-invariant) section;
\item there exists a positive-definite, norm-preserving, antiholomorphic
bundle
map $\tilde\tau:\Cal{V}\to \Cal{V}^*$  lifting $\tau$.\endroster

Conversely, every such $\Cal{V}$ defines a solution of the Bogomolny
equations.\endproclaim
\noindent (Note that this theorem is a generalisation of
\cite{\Hi, Theorem 4.2}, which is stated in a way specific to $\SU
2$.)

We have encoded unitons as solutions $(\D{},\Phi)$ of the
Bogomolny equations with special properties (time invariance, finiteness,
reality).   We will show in \S4 that these properties correspond to
properties of the bundle  $\Cal{V}\to T\P1$.  Namely,
\roster
\item  time translation induces a one-complex-parameter family of
automorphisms of $T\P1$ which lift to bundle maps.
\item Finiteness translates as an extension of the bundle to the fibrewise
compactification of $T\P1$.  The time translation $\delta_{t}$
fixes the divisor at infinity and the bundle lift
$\widetilde\delta_{t}$ is the
trivial bundle map over the divisor at infinity.
\item Reality of the uniton solutions translates as a real structure (a
fixed antiholomorphic principal bundle involution $\tilde{\sigma}$ over
the antiholomorphic involution $\sigma$ of $T\P1$, or a lift
to a map $\Cal{V}\to \Cal{V}^*$) which is positive
definite above a fixed point.
\endroster

\subhead \ardvark\ Inverse Construction\endsubhead
Conversely, given such a bundle we can construct the Bogomolny solution.

{}From the sections map
$$\align \P1\times\C^3&\to\OO{}(2)\cong T\P1\\
\l,(a,b,c)&\mapsto\eta=\frac12a-b\l-\frac12c\l^2\endalign$$
we get a pullback of the bundle $\Cal{V}$ to $\P1\times\C^3$.  Over the
open
set $Y$ of sections over which $\Cal{V}$ is trivial, we can push this
bundle
down from $\P1\times Y$ to a trivial bundle over $Y$.
Call this bundle $\Cal{W}\to Y$.

Now put the quadratic form $(db)^2+(da)(dc)$ on the holomorphic tangent
space to $\C^3$.  To a complex metric we can associate null planes on
which the restricted metric is degenerate, and null lines on which the
restricted metric is zero.  Each null line lies in a unique null plane
(its orthogonal complement).
  The null planes are parametrised by $p\in T\P1$ and
given by the space of sections of $T\P1$ through $p$, i.e.\ if
$p=(\l,\eta)$ then $(a,b,c)$ are constrained by
$\eta=\frac12(a-2b\l-c\l^2)$.  Restricted to this plane, the metric
$(db)^2+(2\l db+\l^2dc)(dc)=(db+\l dc)^2$ is degenerate.  Each null plane
$\Pi_p$ inherits a flat `null' connection:  Let $p\in T\Bbb P^{1}$ be the
unique point of intersection of the family of sections $\Pi_p$.  A
fixed frame of $\Cal{V}|_p$ induces a frame of $\Cal{W}|_{\Pi_p}$.
Define
$\D{\Pi_{p}}$ to be the `null' connection for which this frame is covariant
constant.

Now fix a point $y\in Y$.  Some directions (lines through $y$ in $\C^3$)
may lie in
two different null planes (they correspond to sections intersecting in
two distinct points), but null lines lie in unique null planes, so we can
define a holomorphic connection on the null lines without ambiguity.
This connection in the null lines determines a unique connection along
all directions.  (See \cite{Hi}.)

There is a second way of defining the connection.  On any line in $\C^3$,
pick two points $y_1,y_2$.  They correspond to two sections
of $T\P1$. Since $T\P1$ is the total space of $\OO{\P1} (2)$, they
intersect in two points $p_{1}$, $p_{2}$ (with multiplicity).
We get a double point \Iff
the sections are tangent at a point \Iff the line was null.
To each point $p$ is associated a connection on the corresponding
hyperplane
 $\Pi_{p_{i}}\subset\C^3$.  The original line will be the intersection of
 the two null planes unless it was a null line, in which case it is only
 contained in the intersection.
Taking the average of the two null connections,
we get a connection on the original line, i.e.\
$$\D{\text{line containing
}y_1,y_2}=\frac12(\D{\Pi(p_1)}+\D{\Pi(p_2)}).\tag3.7$$
In the
case of a null line, we get back the null connection, because any two
sections will intersect in a double point.

Since this
construction also gives a holomorphic connection, which agrees
with the first on null lines, it follows from the preceding discussion
that they are identical.

\subhead 3.8 Higgs' Field\endsubhead
As we defined $\D{}$ along a line as
$1/2(\D\Pi+\D{\tau \Pi})$, where the line is the intersection of a
null plane and its conjugate, we can also define
$$\Phi d\chi=\frac{i}{2}(\D\Pi-\D{\tau
\Pi}),\tag3.9$$
where the coordinate $\chi$ on the line of intersection can be well
defined. (See \cite{{Hu1}}.)

\bigpagebreak
In \S\S7-10, we will use complex
algebraic methods
to understand what happens in the finite (energy) case.
The space of lines, $T\P1$, can be embedded into $\P3$ as a
quasi-projective variety.  It can be compactified
by adding a singular point.
Since $\Cal{V}$ is trivial over the section at infinity, we get a
bundle over this variety (which is a degenerate conic).
Sections of $T\P1$ correspond to certain hyperplane sections.
In fact, the hyperplane sections of $\P3$ are given by $(\P3)^*$, and the
sections of $T\P1 \to \Bbb P^{1}$ correspond to the set of hyperplane
sections which do not meet the singular point
$\Bbb C^{3}\subset \Bbb P^{3}$.
We work out that the points we need to add from $\P2=\P3\setminus\C^3$
are just the set of
hyperplane sections restricted to which the bundle is trivial, so we do in
fact get back a (finite) uniton.

\head4. Adapted Coordinates\endhead

In this and the next sections we assume the pair $(\D{},\Phi)$
comes from a uniton in the way specified in \S2.

Left-inverse to stereographic projection is the embedding
$\P1\hookrightarrow\R3$:
$$\lambda\overset i\to\mapsto \left(\frac{\l+\lb}{1+\l\lb},
-i\frac{\l-\lb}{1+\l\lb},\frac{1-\l\lb}{1+\l\lb}\right)\tag4.1$$
alternatively, $\P1\hookrightarrow\C\times\R{}$:
$$\lambda\overset i\to\mapsto
\left(\frac{2\l}{1+\l\lb},\frac{1-\l\lb}{1+\l\lb}\right).$$
Using this inclusion we get an exact sequence of (real) bundles over
$\P1$:
$$0\rightarrow T\P1\overset i_*\to\rightarrow
\left.T\R3\right\vert _\P1\cong\P1\times\R3\rightarrow N_\P1\rightarrow
0,$$
where $N_\P1$ is the normal bundle of the embedding,
and
$$\align
i_*\left(\dd\l\right)
  &=\frac{2}{(1+\l\lb)^2}\dd z
  -2\frac{\lb^2}{(1+\l\lb)^2}\dd \z
  -2\frac{\lb}{(1+\l\lb)^ 2}\dd t \\
i_*\left(\dd\lb\right)
  &=-2\frac{\lb^2}{(1+\l\lb)^ 2}\dd z
  +\frac2{(1+\l\lb)^ 2}\dd \z
  -2\frac\l{(1+\l\lb)^ 2}\dd t.
\endalign $$
The isomorphism $\R3\times\P1\cong T\P1\oplus N_\P1$
suggests that we find adapted coordinates to replace $z$, $\z$, $t$,
$\lambda$ on $\R3\times\Bbb P^{1}$.  The fibre
coordinate of the trivial real line
bundle $T\P1\oplus N_\P1\rightarrow T\P1$
will be $u$, given by
$$\lambda,u\in\Bbb{P}^1\times\Bbb{R}
 \mapsto(\lambda,ui(\lambda))\in\P1\times\R3.$$
On the complement,
$\eta i_*(\dd\l)$ gives us the fibre coordinate of
$T\P1\rightarrow\P1$.
The coordinate systems for a given $\lambda$ are related linearly by
$$\eebu=\pmatrix
1/2 &-\l^2/2 &-\l\\
-\lb^2/2 &1/2&-\lb\\
\lb/(1+\l\lb) &\l/(1+\l\lb) &(1-\l\lb)/(1+\l\lb)
\endpmatrix\zzbt.\tag\brontosaurus$$
Remark that
$\eta=\frac12(z-2\l t-\l^2\z)$,
$\eb=\frac12(\z-2\lb t-\lb^2z)$, and restricted to the plane $\{u=0\}$,
$z=\frac{2\eta-2\l^2\eb}{(1+\l\lb)^2}$,
$\z=\frac{-2\lb^2\eta+2\eb}{(1+\l\lb)^2}$,
$\eta=\frac12\frac{(1+\l\lb)(z+\l^2\z)}{1-\l\lb}$,
$\eb=\frac12\frac{(1+\l\lb)(\lb^2 z+\z)}{1-\l\lb}$
(away from $\{|\l|=1\}$).
{}From this change of coordinates, we can relate  the
connection on $\R3$ in the two coordinate systems.

Although we have not made the distinction,
changing the $\R3$ coordinates in a $\l$ dependent way also affects
$\dd\l$.
To be precise, we should have used coordinates $\l,\lb,z,\z,t$ and
$\l',\lb',\eta,\eb,u$, with $\lp=\l$.
The distinction will be important when we want to show that the
$\dbar$-operator extends to $\eta\in\P1$ in some neighbourhood of $\l=0$,
because we will need to work with $(\eta,\eb,u)$-coordinates.

We make use of the fact that if $\zeta$ and $\chi $ are two choices of
coordinates, $d\zeta=\Lambda d\chi\iff\dd\chi=\Lambda^t\dd\zeta.$
We have calculated $\eebu=B\zzbt$
above.  It follows that
$$d\pmatrix\l'\\ \lbp\\ \eta\\ \eb\\ u\endpmatrix=\pmatrix
1 & 0 & 0 & 0 & 0\\
0 & 1 & 0 & 0 & 0\\
 \\
\dB\l & \dB\lb & & B &\\
{}&{}\\
\endpmatrix d\pmatrix\l\\ \lb\\z\\ \z\\ t\endpmatrix,$$
and for example
$$\align\dd\lb&=\dd\lbp+\(\dB\lb\)^t\pmatrix\dd\eta\\ \dd\eb\\ \dd
u\endpmatrix\\
&=\dd\lbp+\(\frac{2\l}{(1+\l\lb)}\eb-u\)\dd\eb+\frac{2}{(1+\l\lb)^2}\eta
\dd u.\endalign$$
So an element of the kernel of the three operators
$\D u-i\Phi,\ \D\eb,\ \dd\lb$ is also in
the kernel of
$$\D\lbp\deq\dd\lbp-\(\frac{2\l}{(1+\l\lb)}\eb-u\)A_\eb
-\frac{2}{(1+\l\lb)^2}\eta (A_u-i\Phi).$$
In these coordinates, recalling \thetag{\parana}
$$\aligned\D u-i\Phi &=\frac1{1+\l\lb}\left[ 2\l\dd z +2\lb\dd\z
+(1-\l\lb)\dd t+2(1+\l)(A_z+\lb A_\z)\right]\\
\D {\eb} &=\frac2{(1+\l\lb)^2}\left[-\l^2\dd z +\dd\z-\l\dd t
  +(1+\l)(-\l A_z+A_\z)\right]\\
\dd\lb&=\dd\lb\\
\D\lbp&=\dd\lbp-\frac{2(1+\l)}{(1+\l\lb)^3}
\((2\eta-2\l^2\eb+\l(1+\l\lb)u)A_z \right.\\ &\quad\quad\quad\left.
 +(2\lb\eta+2\l\eb-(1+\l\lb)u)A_\z\).\endaligned\tag\frog$$
Sections of the bundle $\Cal{V}$ correspond to elements of the intersection
of the kernels of these
operators.  Roughly speaking, the system has enough solutions if it is
involutive (see \cite{\Wr}).
\proclaim{Lemma 4.4} The system $$\{\D u-i\Phi,\D \eb,\dd\lb\}=\{\D
u-i\Phi,\D
\eb,\D\lbp\}$$
is involutive \Iff $S$ is harmonic.  \endproclaim
\demo{Proof}
A system of differential operators given by generators is involutive \Iff
all Lie brackets of pairs of generating operators lie in the system, i.e.\
in the span of the generators.
Of course any other set of generators is as good, and finding alternative
generators which have vanishing Lie brackets makes things simpler.
That said, remark that
$$\aligned\frac12\{(1+\l\lb)\D\eb+\l(\D u-i\Phi)\}
&=\dd\z-\frac12\l\dd t+(1+\l)A_\z, \text{ and}\\
\frac1{2\l}\{-\lb(1+\l\lb)\D\eb+(\D u-i\Phi)\}
&=\dd z+\frac1{2\l}\dd t+(1+\l^{-1})A_z.\endaligned\tag\aligator$$
It follows that $[\D\eb,\D u-i\Phi]=0$ \Iff
$[\dd \z+(1+\l)A_\z,\dd z+(1+\l^{-1})A_z]=0$.  But since $A_z$ and $A_\z$
do not depend on $t$, this is the case \Iff $S$ is harmonic.  That the
entire
system is involutive follows from the fact that \thetag\aligator\ does not
depend on $\lb$.\qed
\enddemo

In fact, the parametrised system of connections $(\dd \z+(1+\l)A_\z,\dd
z+(1+\l^{-1})A_z)$ having curvature zero can be trivialised over
$\{(\l,z)\in\Cs\times {\Bbb S}^2\}$, giving  Uhlenbeck's  extended
solutions.

\head 5. Compactness\endhead
We are interested in extending the bundles from $T\P1$ to $\T$.
This will follow from the finiteness (of the energy) of the uniton.
The problem is that $\D{}=d+A$, which depends on $(\lambda,\eta,u)\in
T\P1\times \R{}$,
does not have a limit as $\eta\rightarrow\infty$ (in the $1/\eta$
coordinate), so the holomorphic structure must be defined differently.

We do so in two stages,
first for $\{0\ne\l\ne\infty\}$, the `nonpolar' fibres of $\T\to\P1$,
 then in neighbourhoods of the poles (i.e.\
$\l\in\{0,\infty\}$).
The first step is motivated by the geometry of the problem and uses the
extended solution, the second
relies on Sobolev methods to give the existence of a continuous gauge
in which the $\dbar$-operator is smooth.

If $S:\P1\rightarrow \U N$ is a uniton, both $S(x,y)=S(z,\z)$ and
$S(\zh,\zhb)=S(\xh,\yh)$ are continuous,
where $\zh=1/z$ etc.  In terms of $A_z$ this means
$$A_z=-\frac1{z^2} A_\zh,\quad A_\z=-\frac1{\z^2} A_\zhb,\tag\panda$$
i.e.\ $A_z$ has a strong vanishing property as $z\ra\infty$.  Writing this
out in terms of $x$ and $y$,
we see that $A_z$ vanishes to order $2-\epsilon$ at infinity.
In geometric terms, it means that the `energy' of the connection is
concentrated around the $t$-axis in $\R3$ (see \figa), so that when
$\lambda\ne 0,\infty$, solutions to $\Du$ should have limits
as
$u\ra\infty$.
The limit as $u\to\infty$ gives a natural holomorphic framing over
\lq nonpolar\rq\  fibres, $P_{\lambda}$, $\lambda\ne 0,\infty$,
which extends to $\eta=\infty$, giving us the
compactification there.

\apict

Because $\{\D\eb,\D\lbp,\D u-i\Phi\}$ is involutive, solutions to this
system locally correspond to solutions of $\{\D\eb,\D\lbp\}$ restricted to
a plane $\{u=u_0\}$ for some $u_0$.
As the figure suggests, away from the poles this makes sense for
$u_0=\infty$.
Near the poles, however, this doesn't work and we choose $u_0=0$.
The transition between the resulting frames amounts to integrating $\D
u-i\Phi$ from $u=0$ to $u=\infty$.

\subhead 5.2 ${G}_\infty$ trivialisation\endsubhead
This geometric line of reasoning can be made precise analytically (see
\cite{\thesis}),
but as we said, the extended solution encodes the analytic information of
the uniton, and can be exploited here. In fact, pursuing the geometric
argument we have sketched uncovers the role of
the extended solution. Now we restrict our attention to
$\lambda\in\Bbb C^{*}$, i.e.\ to nonpolar fibres.

We saw that $\{\D \eb,\, \D u-i\Phi\}$ has the same kernel as
\thetag{\aligator},
$$\align &\dd\z-\frac12\l\dd t+(1+\l)A_\z\\
         &\dd z+\frac1{2\l}\dd t+(1+\l^{-1})A_z\endalign$$
away from the fibres ($\l=0,\l=\infty$) where \thetag{\aligator} has
singularities.  Both systems are \lq underdetermined\rq, since any
solution can be multiplied by any holomorphic function in $\eta$ to give
another solution.

The exact function, however, may be fixed by adding another differential
operator to the system.
Over $\{\l\in\Cs\}$, we can add $\dd t$ and get a completely integrable
system on any fibre $\{\l=\l_0\}\subset\T\times\R{}$, $\l_0\in\Cs$, i.e.\
a (full) smooth
connection on $\Cal{V}$ restricted to a fibre of
$\T|_\Cs\times\R{}\ra\Cs$, with
zero curvature.

Since solutions of the augmented system are independent of $t$, we can
push the system down to $\{(z,\z)\in\R2\}$; i.e.\ all solutions to the
augmented system are obtained by pulling back solutions of
$$(\dd z+(1+\l^{-1})A_z,\dd\z+(1+\l)A_\z)\tag\lobster$$
on $\R2$ which are none other than the extended solutions.

Since the lines $\{u\in \Bbb R, \eta, \lambda\ \text{ constant}\}$
project to
circles on ${\Bbb S}^2$ (\figg), and the extended solution is defined on
${\Bbb S}^2  \times \Bbb C^{*}$, the kernel of $(\nabla_{\bar{\eta}},
\nabla_{\bar{\lambda^{\prime}}}, \nabla_{u} - i\Phi)$ extends to
$$\{(\lambda, \eta, u) \in \Bbb C^{*} \times
  \Bbb C \times \Bbb S^{1}\}.$$
  If we
require our extended solutions to be based at $\infty$, i.e.\
$E_{\lambda}(\infty) = \one$, then the lifted solution is
constant $(\one)$ at $u=\infty$, and we can take this to define the
${G}_{\infty}$ trivialisation of $\Cal{V}$ restricted to nonpolar fibers
of $\T$.

Similarly, away from $\{|\lambda| = 1\}$, the equator, $\eta$-planes
$$\{\eta\in\Bbb C,\lambda = \lambda_{0}, u = u_0\}$$ are mapped to
$\Bbb S^{2} \setminus
\{\infty\}$, so the lifted extended solution extends to $$\{(\lambda,
\eta, u) \in \Bbb C^{*} \times \Bbb P^{1} \times \Bbb S^{1}: |\lambda|
\ne 1\}.$$  (Concentric circles $\{|\eta |=r,\lambda =\lambda_{0},
u=u_{0}\}$ are mapped to nested circles on $\Bbb S^{2}$ containing
$\infty$.)
Parallel translation by $\nabla_{u} -i \Phi$ from $\{u=\infty\}$ to
$\{u=u_{0}, \text{finite}\}$ will give the transition between the
${G}_{\infty}$
trivialisation of the nonpolar fibres and the extension of the bundle in
polar neighbourhoods which we define presently.

\gpict

\subhead 5.4 Polar neighbourhoods \endsubhead
In \S3 we defined the holomorphic structure of $\Cal{V}\mapsto\T$ by
specifying
a $\bar{\partial}$-operator.  We have just seen that for
$\lambda\in\Bbb C^{*}$, the extended solution integrates this
operator to give a trivialisation of
$\Cal{V}|_{\widetilde{T \Bbb C}^{*}}$.
Over \lq polar\rq\
neighbourhoods ($\lambda$ near $0$ or $\infty$)
we show that the $\bar{\partial}$-operator extends to a
smooth operator near $\eta=\infty$
(equivalently $\hat{\eta} = \infty$), after a
continuous change of gauge, showing that bundles associated to unitons
extend to the fibrewise compactification---as holomorphic bundles. (A
priori we don't know that the $\bar{\partial}$ is integrable near
${G}_{\infty}$; this follows from smoothness.)
The sets
$$\aligned
U_0 & =\{(\l,\eta):\l\ne\infty,|\eta|\leq\infty\}\\
 U_\infty & =\{(\l,\eta):|\l|<1/2,\eta\ne0\}\\
\hat U_0 & =\{(\lh,\eh):\lh\ne\infty,|\eh|\leq\infty\}\\
\hat U_\infty & =\{(\lh,\eh):|\lh|<1/2,\eh\ne0\}\\
U & =\{(\l,\eta):0\ne\l\ne\infty\}=\Cs\times\P1,\endaligned\tag5.5$$
form a cover of $\T$.  Let
$$\ue=\{|\ebp|\le1\},\quad
\ul=\{|\l|\le1\}.\tag5.6$$

We will smooth the $\bar{\partial}$-operator on $U_{\infty}$.  The
situation on $\hat{U}_{\infty}$ is similar.  On $U_{\infty}$ the
$\bar{\partial}$-operator has components
$$\align
\D\ebp &=-\eb^2\D\eb=\dd\ebp-\frac{(1+\l)(1+\l\lb)^2}2 \{\frac\l{(\l^2-
\ebp/\ep)^2}A_\zh-
\frac1{(1-\lb^2\ebp/\ep)^2}A_\zhb\}\\
&\deq\dd\ebp+B_0A_\zh+B_1A_\zhb,\tag5.7a\\
\D\lbp&=\dd\lbp-\frac{2(1+\l)}{(1+\l\lb)^3}
  \[\frac\ebp{1-\l^2\ep/\ebp}A_\zh
  +\frac{\lb\ebp{}^2/\ep+\l\ebp}{(1-\lb^2\ebp/\ep)^2}A_\zhb\],\tag5.7b
  \endalign$$
where for simplicity we restrict to the hypersurface $u=0$.

Both operators are smooth away from $\eta'=0$, and near $\eta'=0$ are
smooth functions of $\eta',\ebp$ and $\ebp/\ep$.  Since $A_\zh$ and
$A_\zhb$ are smooth near $\ep=0$ (for $\l\in\ul$), the degree of
smoothness of the $\bar{\partial}$-operator is determined by that of the
coefficients $B_i$.  For example, the $(\ebp/\ep)$ factor is bounded and
hence locally lies in $(L^{2})$, but its derivative does not.

\smallpagebreak

Atiyah and Bott have determined when local smoothing gauges exist
in one complex dimension
(\cite{\AB, 5.1,\S14}).
By approaching the smoothing as a
parametrised one-dimensional problem and taking advantage of the special
form of the singularity, we will be able to find a continuous holomorphic
gauge.  Such gauges might not exist in general for $\dbar$-operators with
coefficients in $L^{2}_{0}$ but without this particular
singularity structure.

\medpagebreak Because all the objects we will be dealing with, eg.\ $B_i$,
are smooth away from $\ep=0$ and depend smoothly on $\lambda\in\Lambda$,
integrability ($L^p_k$, i.e.\ $L^p$
integrability of partial derivatives up to order $k$) on $\ue\times\ul$,
and
on fibres of $\ue\times\ul\to\ul$ are equivalent.  In fact, $B_i\in
L^2_0(\ul\times\ue,\gl{N})$ can also be seen as a smooth map valued
in a function space: $$\Cal B_i\in C^\infty(\ul,L^2_0(\ue,\gl{N}
).\tag5.8$$ By
taking the second view of $B_i$, as a smooth function valued in a function
space, we reinterpret the search for a smooth gauge as a parametrised
one-complex-dimensional problem.

The basic tool for proving smoothness is the \proclaim{Sobolev Lemma}
There
are inclusions
$$\align L^2_2(\ue,\gl N)&\subset C^0(\ue,\gl N)\text{ and}\\
  L^2_3(\ul\times\ue,\gl N)&\subset C^0(\ul\times\ue,\gl N) \endalign$$
which are continuous with respect to the Sobolev and supremum norms
respectively.\endproclaim For our purposes, continuity of the inclusions
will be very important.  See \cite{\GH,\,p86} for a proof.

\proclaim{Lemma 5.9} The operator $P:g\mapsto\dd\ebp g g^{-1}$ can be
extended
to a smooth invertible map $$P:\L 2 k(\ue,\GLC N)_0\ra\L 2{k-1}(\ue,\glc
N)$$
for $k>2$, where $\L 2 k(\ )_0$, indicates the space of based maps,
$g(0)=\one$.\endproclaim

\demo{Proof} $P$ extends to a map of Sobolev spaces
because \roster\item since $\ue$ is compact and $L^2_k(\ue)\subset
C^1(\ue)$, we can find a constant such that $\norm
{g^{\prime-1}}<\text{const}\norm g'$ for $g'$ in some neighbourhood of
$g_{0};$
\item $L^p_k(\R n)$ is a Banach algebra for $k>n/p$ and $L^p_j$ is a
topological $L^p_k$- module for $0\le j\le k$ \cite{\AB,14.5};  \item
$\dd\ebp$ gives a map $\L 2 k(\C^n,\GLC N)\ra\L 2{k-1} (\C^n,\glc N)$ for
all $k,n$.\endroster

We can calculate the derivative, $DP$, of $P$ by expanding
$$\align P\(g_0(\one+g_1)\)&=\dd\ebp g_0g_0^{-1}+\dd\ebp
g_1g_0^{-1}-g_1g_0^{-1}\dd\ebp g_0g_0^{-1}+\phi_{g_0}(g_1),\\
\intertext{where $\phi_{g_0}(g_1)$ is tangent to the zero map (i.e.\
$\lim_{|g_1|\ra0}\frac{\phi_{g_0}(g_1)}{|g_1|}=0$),}
&=P(g_0)+DP(g_0)(g_1)+\phi_{g_0}(g_1).\endalign$$
In particular $DP(\one)=\dd\ebp$, so we can apply the inverse function
theorem for Banach spaces \cite{\La, I.5.1} to get an inverse to $P$ in a
neighbourhood of $P(\one)=0$.
It is smooth because $P$ is smooth:

We can verify the existence of higher derivatives for $P$ either by
iteratively differentiating $P$, or by using the chain rule and remarking
(a) that $\dd\ebp$ and $m(a,b)=ab$ are linear and multilinear respectively
and hence both smooth \cite{\La, I.3.12}; and (b) that $g\mapsto g^{-1}$
is smooth because its $k\orth$ derivative at $g=g_0$ is the $k$-linear
map
$\bigoplus^k L^2_k(\ue,\glc N)\ra L^2_k(\ue,\glc N)$, given by

$$(g_1,...,g_k)\mapsto\sum_{\sigma\in
S_k}g_0^{-1}g_{\sigma(1)}g_0^{-1}g_{\sigma(2)}g_0^{-1}\cdots
g_{\sigma(k)}g_0^{-1}.$$

Note that we have made repeated use of (1) and (2) above.
\qed\enddemo

Looking for a better gauge involves integrating $A_\ebp$, i.e.\
finding $C_0$
and
$C_1$ such that $\dd\ebp C_i=-B_i$,
$$\align
C_0&=\frac{(1+\l)(1+\l\lb)^2}2\frac{-\l\ep}{\l^2-\ebp/\ep},\\
C_1&=-\frac{(1+\l)(1+\l\lb)^2}2\frac\ebp{1-\lb^2\ebp/\ep},\tag5.10
\endalign$$
which are not only continuous at $\ep=0$, but vanish there. Use
$$g=\id+C_0A_\zh+C_1A_\zbh,$$ to change gauge
$$A_\ebp^g\deq g^{-1}A_\ebp
g+g^{-1}\dd\ebp g$$ giving $A^{g}_{\bar{\eta}^{\prime}}$ the same
continuity
properties :
$$\align
A^{g}_{\bar{\eta}^{\prime}} = (B_0 A_\zh+B_1 A_\zbh)^g&= g^{-1}(B_0
A_\zh+B_1 A_\zbh)g+g^{-1}\dd\ebp g\\
&=B_0 A_\zh+B_1
  A_\zbh+C_1B_0[A_\zh,A_\zbh]+C_0B_1[A_\zbh,A_\zh]\\
&\quad-\(C_0A_\zh+C_1A_\zbh\)\(B_0 A_\zh+B_1
  A_\zbh\)\(C_0A_\zh+C_1A_\zbh\)\\
&-B_0A_\zh-B_1A_\zbh
  +C_0\(B_2\dd\zh A_\zh+B_3\dd\zhb A_\zh\)\\
&\quad+C_1\(B_2\dd\zh A_\zbh+B_3\dd\zhb A_\zbh\)
  -\(C_0A_\zh+C_1A_\zbh\)\\
&\quad\quad\[-B_0A_\zh- B_1A_\zbh+C_0\(B_2\dd\zh A_\zh+B_3\dd\zhb
  A_\zh\)\right.\\&\quad\quad\quad\left.
  +C_1\(B_2\dd\zh A_\zbh+B_3\dd\zhb A_\zbh\)\],
\endalign$$
where the terms $B_2,B_3$ arise because $\dd\ebp\ne\dd\zbh$:  (on
$u=0$)
$$\align\dd\ebp&=-\eb^2\dd\eb\equiv-
\frac{\l^2(1+\l\lb)^2}{2(\l^2-
\ebp/\ep)^2}\dd\zh+\frac{(1+\l\lb)^2}{2(1-
\lb^2\ebp/\ep)^2}\dd\zbh\pmod{\dd t}\\
&\deq B_2\dd\zh+B_3\dd\zbh\pmod{\dd t}.\endalign$$
We can ignore the $\dd t$ terms because the connection coefficients are
independent of $t$.

Since each (persistant) term contains a factor $C_j$, $\Aeg$ vanishes on
$\ep=0$, and is continuous.
Because $\dd\ebp C_j=B_j$, $\dd\ebp\Aeg$ is bounded but discontinuous, we
see that $\Aeg\in\L21(U_\infty,\glc N)$, in fact
\proclaim{Lemma 5.11}  The map
$$\align \CalAeg:&\ul\ra\L21(\ue,\glc N)\\
\intertext{such that}
   & \CalAeg(\l)(\ep_0)=\Aeg(\l,\ep_0)\endalign$$
is smooth.\endproclaim
\demo{Proof}
We have to show (1) that the function $\Aeg|_{\l=\l_0}$ and its first
$\ep,\ebp$ derivatives are square
integrable, for all $\l_0\in\ul$, and (2) that $\dd\l^p\dd\lb^q\CalAeg$
exist and are in $L^2_1$ for all $p$ and $q$.

(1) Since
$$\Aeg=g^{-1}A_\ebp g+g^{-1}\dd\ebp g$$
is smooth away from $\ep=0$, and has a singularity of type
$\ep\phi(\ep/\ebp)$ there, its first derivatives in $\ep$ or
$\ebp$ may have a bounded discontinuity,
which doesn't affect the finiteness of the $L^2$ norm.
In fact, we can multiply $\Aeg$ by the complement of a bump function of
arbitrarily small mass concentrated at $\ep=0$ and find that the map
$\CalAeg:\ul\ra\L21(\ue,\glc N)$ is continuous.

(2) Since $\Aeg$ is smooth on $\{\ep\ne0\}$ we can take the $\l$
derivatives
of $\CalAeg$ pointwise, i.e.\ when $\ep\ne0$
$$\dd\lp^p\dd\lbp^q\CalAeg(\l_0)(\ep_0)
=\dd\lp^p\dd\lbp^q\Aeg|_{(\l_0,\ep_0)}.$$
Examining the terms of $\Aeg$ one at a time, we find that all partials
are continuous.
For example,
$$\align \dd\lp
C_0&=\frac{(1+\l\lb)^2+2\lb(1+\l)(1+\l\lb)}2\frac{-\l\ep}{\l^2-\ebp/\ep}\\
  &\quad + \frac{(1+\l)(1+\l\lb)^2}2
    \frac{-\ep(\l^2-\ebp/\ep) +2\l^2\ep}{(\l^2-\ebp/\ep)^2},\\
\dd\lp C_1&=
\frac{(1+\l\lb)^2+2\lb(1+\l)(1+\l\lb)}2
\frac\ep{\lb^2(1-\lb^2\ebp/\ep)}.\endalign$$

Since all its partials exist and are continuous,
$\CalAeg:\ul\ra\L21(\ue,\glc
N)$ is smooth.\qed\enddemo

Now we can exploit the fact that $P$ has a {\it smooth} inverse.
As a result, $$\gt=P^{-1}\circ\CalAeg:\ul\ra\L22(\ue,\GLC
N)$$ is a smooth map and $P(\gt)=\Aeg$.
Composing with the continuous Sobolev embedding $L^2_2\hookrightarrow
C^0$, we see that $\gt$ is a {\it continuous\/} change of gauge over
$U_\infty$, such that the $\ebp$-operator in this gauge is trivial,
i.e.\ $\Aegt=0$, and $\Algt=\gt^{-
1}\Alg\gt+\gt^{-1}\dd\lb\gt$ is continuous (all its ingredients are).

\bigpagebreak
Since $0=[\D\lb,\D\ebp]=[\dd\lb+\Algt,\dd\ebp]$,
it follows that $\Algt$ is
 meromorphic in $\ep$.  Since it is continuous, it is holomorphic.
Using the fact that $\Algt$ is smooth near $\{|\ep|=1\}$ and
differentiating the Cauchy integral
$${\dd\ep}^j{\dd\l}^k{\dd\lb}^l\Algt(\l_0,\ep_0)
=\int_{|\ep|=1}\frac{\left({\dd\l}^k{\dd\lb}^l\Algt(\l,\ep)\right|_{\l
=\l_0}}{(\ep-\ep_0)^j},$$
we see that $\Algt$ is smooth on
$\{|\ep|<1\}$.
We can then find a smooth change of gauge $\gh$ such that
$\Algh=0=A_\ebp^{g\gt\gh}$.
Then, $s_\infty\deq g\gt\gh$ is the required holomorphic trivialisation
over
$U_\infty$.

A similar construction works over the south pole, completing the proof
that bundles arising from unitons extend to $\T$.

\head 6. Extra Structure\endhead

\subhead 6.1 Triviality over the section at infinity\endsubhead
Lifting a based extended solution to the total space $(\Bbb R^{3} \times
\Bbb P^{1})$ gives a section $\tilde{f}$, of $\Cal{W}|_{\Bbb R^{3} \times
\Bbb C^{*}}$, which pushes down to a choice of gauge on
$\Cal{V}|_{\widetilde{T
\Bbb C}^{*}}$.  Let $f_{\infty}$ and $f_{\widehat{\infty}}$ be holomorphic
gauges of $\Cal{V}$ over $U_{\infty}$ and $\hat{U}_{\infty}$ respectively
\thetag{\lobster}, and
let $\tilde{f}_{\infty}$ and $\tilde{f}_{\widehat{\infty}}$ be their
pullbacks to corresponding regions of the total space.  Since
$\lim_{\eta\mapsto\infty} \nabla_{\bar{\lambda}^{\prime}} =
\dd{\bar{\lambda}^{\prime}}$ on every plane $\{u=\text{constant}\}$,
$\tilde{f}_{\infty}|_{\{\eta^{\prime} = 0\}}$ and
$\tilde{f}_{\widehat{\infty}}|_{\{\hat{\eta}^{\prime} = 0\}}$
are holomorphic
functions in $\lambda$.  Without loss of generality we can assume
$\tilde{f}_{\infty}|_{\{\eta^{\prime} = 0\}} \equiv \one \equiv
\tilde{f}_{\widehat{\infty}}|_{\{\hat{\eta}^{\prime} = 0\}}$.
In terms of the
chosen gauges, $\Cal{V}$ has transition matrices
$\tilde{f}^{-1}\tilde{f}_{\infty}$ and
$\tilde{f}_{\widehat{\infty}}\tilde{f}$ (which are
independent of $u$ because
$\tilde{f}$, $\tilde{f}_{\infty}$ and $\tilde{f}_{\widehat{\infty}}$ all
solve $\nabla_{u} -i \Phi$).  To calculate
$\tilde{f}|_{\{\eta = \infty\}}$, we
appeal to the map
$\Bbb R^{3} \times \Bbb C^{*} \mapsto {\Bbb S}^2 \times \Bbb
C^{*}$.
See \figg.

For $\lambda_{0}\in \Bbb C^{*}$, we saw that the lines $\{\eta =
\eta_{0}, \lambda = \lambda_{0}\}$ (i.e.\ $u \in \Bbb R$) are mapped
to circles on ${\Bbb S}^2$ through  $\infty$, with $\infty$ missing.
Similarly, for $|\lambda_{0}| \ne 1$, the planes $\{u=u_{0}, \lambda
=\lambda_{0}\}$ are mapped to ${\Bbb S}^2\setminus\{\infty\}$.
Since $\tilde{f}$
was pulled back from the based extended solution $E_{\lambda}$ on $\Bbb
{\Bbb S}^2 \times \Bbb C^{*}$ with $E_{\lambda}(\infty)=\one$,
$\lim_{\eta\mapsto\infty} \tilde{f} = \one$ (away from $|\lambda| =1$).
Hence $\Cal{V}|_{{G}_{\infty}}$ has constant transition matrices in these
gauges
and must be trivial.

\subhead 6.2 Time Invariance\endsubhead
Time translation $(z_0,t_0)\mapsto(z_0,t_0+ t)$ induces a one-parameter
group of transformations of $T\P1$.
In coordinates, $(\l,\eta)\mapsto(\l,\eta-t\l)$.
The coefficients of $A_z,A_\z$ and hence $(\D{},\Phi)$ are
independent of $t$,
i.e.\ they are invariant under the group of
translations of $t$.
So the space of solutions to $\D u-i\Phi,\D\eb,\D\lbp$ is invariant under
time translation.

On $T\P1$, the space of oriented geodesics in $\R3$, time translation acts
by $(\l,\eta)\overset\delta_t\to\mapsto(\l,\eta-t\l)$.
The geodesic itself is shifted with respect to the geodesic parameter $u$,
$$u\overset\delta_t\to\mapsto u+\frac{1-\l\lb}{1+\l\lb}t.$$
So a solution $\st(\l,\eta,u)$ such that $(\D
u-i\Phi)\st=0=\D\eb\st=\D\lbp\st$ generates a family of solutions
$$\st_t(\l,\eta,u)=\st(\l,\eta+\l t,u-\frac{1-\l\lb}{1+\l\lb}t).$$
And the map $\tilde\delta_t:\st\mapsto\st_t$, is a bundle isomorphism
lifting $\delta_t$.

Since $$\lim_{z\to\infty}A=0=\lim_{z\to\infty}\Phi,$$ the bundle map is
just
the identity over the section at infinity.
One can also see this by remarking that the ${G}_{\infty}$ trivialisation
is preserved by the time translation map since its lift to
$\Bbb R^{3} \times \Bbb C^{*}$ is independent of $t$.
We will see that the specification of this map encodes the
time independence of the uniton.

\subhead 6.3 The Real Structure\endsubhead
As remarked in \S5, in adapting Hitchin's construction, there
is some choice as to the real structure.
On $\C^3$ one thinks of the real structure literally as a real slice:  a
three dimensional subspace of the {\it real} six dimensional $\C^3$ which
as a set spans the {\it complex} three dimensional $\C^3$.
Any such set is the fixed set of an  antilinear
involution---the real structure.
The $\R3$ of Hitchin's original construction
is the standard real slice
of standard $\C^3$ (with conjugation as the real structure).
Conjugation, however, is not the appropriate real structure for our
purposes.
But $\R3$ is still an {\it invariant} set of the appropriate real
structure:  $$\sigma^{*}(x,y,t)=(\bar x,\bar y,-\bar t=i^{-1}\bar{it}).$$
So when our real structure acts on $\Cal{V}\ra\R3$ it not only conjugates
the
fibres, but reflects $\R3$ in the $x-y$ plane.

One way to understand how the real structure on $\Cal{V}$ arises is to
work with
frames rather than sections.  This is because the real structure on
$\Cal{V}$
comes from the real structure on the complex group ($\Gl N$ in our case),
which induces a real structure on the trivial principal bundle of frames
of
$\Cal{W}$ over $\R3$.  The real structure fixes a real subgroup, and is
$X\mapsto(X^*)^{-1}$ in the case of $\U N$, which is both an involution
and
antiholomorphic with respect to the natural complex structure of $\Gl N$.
A frame of $\Cal{V}$, either locally or at a point, is an invertible
solution,
$\tilde f$, to $$(\D{a\dd x+b\dd y+c\dd t}-i\Phi)\tilde f=0.\tag\fly$$
Everything here lives in $\gl N$, so we can apply the transformation
$$X\mapsto-\overline{\tilde f^{-1}X\tilde f^{-1}}^t$$ to \thetag\fly\ to
get a new equation, which since $iA_t, i\Phi, A_x,A_y\in\u N$, gives
$$(a\dd x+b\dd y+c\dd t+aA_x+bA_y-cA_t-i\Phi)(\tilde f^*)^{-1}=0.$$
Pulling
back by $\sigma$, and using the fact that $A_x,A_y,A_t,\Phi$ are
independent of $t$, we get
$$(\D{a\dd x+b\dd y-c\dd t}-i\Phi)(\sigma^*\tilde f)^{*-1}=0.$$
So the real structure on the
principal bundle $\R3\times\Gl N$ induces an antiholomorphic involution of
the principal bundle of frames of $\Cal{V}\ra T\P1$ which covers $\sigma$
and
which conjugates the natural frames above real sections (i.e.\
(frame)$\mapsto$(frame)${}^{*-1}$ in a unitary frame).  The specific form
of this conjugation is important---other conjugations correspond to
different real groups, i.e.\ $\operatorname{Gl}(N,\R{}),
\operatorname{U}(n,N-n)$  etc.

On the level of vector bundles, this translates into an
antiholomorphic map $\Cal{V}\to \Cal{V}^*$ lifting
$\sigma$.

\subhead 6.5 Framing \endsubhead
Finally, let $p$ be any fixed point of $\sigma$ contained in $P_{-1}$.
We take as the framing above this point the solution to $(\D
u-i\Phi)\phi=0$ along the line $p\subset\R3$ with
$$\lim_{u\to\infty}\phi=\one.$$
Since $p$ is fixed by $\sigma$, the real structure $X\mapsto X^{*-1}$
takes solutions of $\D u-i\Phi$ along $p$ to solutions, and maps our
particular solution to itself.  Since $\Cal{V}|_{P_{-1}}$ is trivial,
$$H^0(P_{-1},Fr(\Cal{V}))@>\text{eval}>>Fr(E_p)$$ is an isomorphism.
This defines the `unitary' framing
$$\phi\in H^0(P_{-1},Fr(\Cal{V}))$$ of the definition.
\remark {Remark 6.6} The value of this result is that since $\T$ is
compact, the a priori analytic uniton bundles are algebraic (see
\cite{{GAGA}}).  This underlies the inverse construction to follow.
\endremark

\head 7. $\T$ as a conic\endhead{}

As a preliminary step to inverting the bundle construction of \S6, we show
that $T\P1$ can be embedded into $\P3$ as
the nonsingular subset of a conic.
Consider the conic $Q$ given in homogeneous coordinates
$\alpha,\beta,\gamma,\delta$ on $\P3$ by $\beta^2=-4\alpha\gamma$.
This conic has a singular point at $[0,0,0,1]$.
Now consider the map $f:T\P1\rightarrow Q$ given by
$$\align
(\l,\eta)&\mapsto[1,-2\l,-\l^2,-2\eta]=[\alpha,\beta,\gamma,\delta]
\tag\antelope\\
(\lh,\eh)&\mapsto[-\lh^2,2\lh,1,2\eh].\endalign$$
The map $f$ extends to a rational map on $\T$ mapping the section at
infinity to the singular point.
Since the bundle is trivial over the section at infinity, when we collapse
this section the bundle descends to another bundle $f_*\Cal{V}$ on $Q$.
More precisely,
\proclaim{Lemma 7.2} Pullback of bundles ($\Cal{Z}\mapsto f^*\Cal{Z}$)
from $Q$ back to
$\T$ is an isomorphism onto the space of bundles on $\T$, trivial over the
section at infinity,
${G}_\infty$.\endproclaim
\demo{Proof}
Pullback of bundles is injective.  (Pushforward is a left inverse of
pullback.)  We  must show it is
surjective---that every bundle trivial on ${G}_\infty$ is the pullback of
a
bundle on $Q$.  Let $\Cal{Z}'\to\T$ be trivial on ${G}_\infty$.  Away
from
${G}_\infty$, $f$ is bijective, so $\Cal{Z}'$ pushes forward to a bundle
on $Q$
away from the singular point.
We shall use the Theorem on Formal Functions to push forward a
trivialisation of $\Cal{Z}'$ in a neighbourhood of ${G}_\infty$ to a
trivialisation of $f_*\Cal{Z}'$ in a neighbourhood of the singular point
$f_*({G}_\infty)$.  So $f_*\Cal{Z}'$ is a bundle (a {\it locally
trivial\/}
sheaf) whose image is $\Cal{Z}'$, proving surjectivity.

 \jpict

Locally, the section at infinity, ${G}_\infty$, looks like the zero
section
of $\O{\P1}(-2)$.
Given local coordinates $(\l,\ep)$ and $(\lh=1/\l,\ehp=\l^2\ep)$ on
$\O{}(-2)$,
a transition matrix for $\Cal{Z}'$ is given over the intersection,
$\{\l\in\Cs\}$, as $$\one+\ep(\phi(\l,\lh,\ep,\ehp)),$$ where $\phi$ is a
polynomial matrix.
Since $\ehp=\l^2\ep$, we can express this in terms of two polynomials as
$$\align\one+&\ep(\phi'(\l,\ep)+\phi''(\lh,\ehp))\\
&=\one+\ep\phi'(\l,\ep)+\ehp\lh^2\phi''(\lh,\ehp),\endalign$$
but not uniquely, as $\ep=\lh^2\ehp$ etc. We can use this property to show
inductively that the bundle must be trivial on all formal neighbourhoods
of
${G}_\infty$, by showing that, for any $k>0$, such a transition matrix
$\in C^1({G}_\infty^{(k)},\GL N)$ is actually a coboundary,
splitting as a product of holomorphic changes of gauge,
i.e.\ that it is in
the
image of $C^0({G}_\infty^{(k)},\GL N)$.

A bundle is trivial on the $(k-1)$st formal neighbourhood,
${G}_\infty^{(k-1)}$, \Iff its transition matrix has the form
$$\one+\ep{}^k(\phi)$$
in some gauge.  Using the fact that $\phi$ can be split as
$\phi=\phi'(\l,\ep)+\phi''(\lh,\ehp)$, we can make a change of gauge:
$$\align(\one-\ep{}^k\phi')&(\one+\ep{}^k(\phi'+\phi''))
  (\one-\ehp{}^k\lh^{2k}\phi'')\\
&=\one+\ep{}^{2k}(\phi'\phi''-\phi'(\phi'+\phi'')-\phi''(\phi'+\phi'')),
\endalign$$
showing that it is trivial on ${G}_\infty^{(2k-1)}$.
Inductively, we get a trivialisation of $\Cal{Z}^{\prime(k)}$, which
is
the same as a maximal rank
section of $\HOM(\E^{\oplus N},\Cal{Z}')$ over the $k\orth$ formal
neighbourhood (${G^{(k)}_\infty}$), for $k$ arbitrarily large,
where $\Cal{E}$ is the trivial bundle.
Now the Theorem on Formal Functions \cite{\Hartshorne, III.11.1} says that
$$f_*\HOM(\E^{\oplus N},\Cal{Z}')^{\wedge}_{[0,0,0,1]}@>\cong>>\varinjlim
H^0(G^{(k)}_\infty,\HOM(\E^{\oplus N},\Cal{Z}')).$$
We have shown that the RHS has a maximal rank element.
The LHS is the set of all sections of $\HOM(\E^{\oplus N},\Cal{Z}')$ on a
neighbourhood of $f^{-1}([0,0,0,1])={G}_\infty$, up to formal equivalence,
i.e.\ germs of sections.
It must contain an element corresponding to the maximal-rank element of
the RHS.  That element is a section on some neighbourhood
of the section at infinity which is nondegenerate on the infinity
section.
Since the determinant function is continuous, it must be nondegenerate on
some neighbourhood where it gives a trivialisation.
This trivialisation pushes forward to give a trivialisation of
$f_*\Cal{Z}'$ in
a neighbourhood of the singular point, so in particular $f_*\Cal{Z}'$ is
a
bundle.\qed\enddemo

The same construction would work with any bundle trivial over a rational
curve of negative self-intersection embedded in a surface, because the
splitting of $\phi$ would go through.  A theorem of Castelnuovo tells us
that the new surface will be smooth \Iff the curve has self-intersection
$-1$, in which case we are just blowing down.

\medpagebreak
The point of this construction is that it tells us what happens to
deformed   real sections in the limit (as $t\toinfty$, for example).
A hyperplane in $\P3$, $\alpha a+\beta b+\gamma c+\delta d=0$, restricted
to $Q$ can be written
$$a\alpha-2\l b\alpha-\l^2\alpha c-2\eta\alpha d=0
\quad \text{or}\quad
-\lh^2\gamma a+2\lh\gamma b+\gamma c z \eh\gamma d=0.$$
When $d\ne0$, we can use affine coordinates $a/d,b/d,c/d$ or just
restrict to the plane $d=1$ and to pull  these sections back to $T\P1$ we
restrict to the affine plane $\{\alpha=1\}$ (see \thetag\antelope).  In
$T\P1$ coordinates, the hyperplane section $(a,b,c,1)$ is
$$\eta=\frac12(a-2b\l-c\l^2),$$
a holomorphic section of $T\P1\to\Bbb P^{1}$.  In other words the
sections $\Bbb P^{1}\to T\Bbb P^{1}$ have their images cut out by
hyperplane sections with $d\ne 0$.

But what about the  $\P2$ of hyperplanes with $d=0$?
{}From \figj, we  see that these are just the
hyperplanes which include the singular point.
Such intersections solve $a\alpha+b\beta+c\gamma=0$ and
$-4\alpha\gamma=\beta^2$, so they solve
$a^2\alpha^2+(2ac+4b^2)\alpha\gamma+c^2\gamma^2=0$.
When $b^2+ac=0$ or $b=0$ the solution is a double line.  In general we get
two lines intersecting in the pinch point:
$$\aligned a^2&\alpha+\left(ac+2b^2\pm2b\sqrt{ac+b^2}\right)\gamma=0,\\
a^2&b\beta=\left(-a^2c+ac+2b^2\pm2b\sqrt{ac+b^2}\right)\gamma.
\endaligned$$
So the correct way to complete the set of holomorphic sections of $T\P1$
is not by adding a section at infinity but by adding a
$\P2$ worth of closed
subvarieties of $\T$, given by the union of the section at infinity and
 two fibres with multiplicity.

\head 8. Compact Twistor Fibration\endhead{}
Let $X\subset \P3\times\P3^*$ be the variety cut out by
$\beta^2+4\alpha\gamma=0$ and $a\alpha+b\beta+c\gamma+d\delta=0$,
where $a,b,c,d$ are homogeneous coordinates on the space of hyperplane
sections of $\P3$, $\P3^*\cong\P3$.
We will use the double (twistor) fibration
$$\matrix
 & & X & & &\\
& \pi_1\swarrow & & \searrow  \pi_2 & &\\
T\P1\hookrightarrow Q & & & & \P3^* \hookleftarrow\C^3.
\endmatrix\tag8.1$$
Pull back the bundle $\Cal{V}$ to $X$ and push it forward to a sheaf
$\Cal W$ over $\Bbb P^{3}$.  Let
$$Y=\{y\in\P3^*:\Cal{V}|_{X_y} \text{ is
trivial}\}\tag8.2$$ and
$$X^{\prime}=\pi^{-1}_{2}(Y).$$
In the following, $y$ will be assumed to be in $Y$.

Grauert's Theorem (\cite{{Gr}}) implies that $\Cal W\to Y$ is locally
free, i.e.\ a bundle, when (the sheaf of sections of) $\pi_{1}^{*}\Cal V$
is flat over $Y$.  Since (sections of) bundles are coherent, this follows
from the flatness of $(X^{\prime}, \Cal O_{X^{\prime}})\to (Y,\Cal
O_{Y})$, which is flat \Iff
$$P_y(m)=\dim_{\C {}} H^0(X_y,\O{X_y}(m))$$ is
independent
of $y\in Y$ for $m\gg0$ (\cite{\Hartshorne, III.9.2}).
We can compute $P_y(m)$ from the long exact homology sequence associated
to the embedding $X_y\subset\P2$:
$$0\ra \O{\P2}(-X_y)\ra\O{\P2}\ra\O{X_y}\ra0.\tag\bandar$$
Since $Q$ is cut out by a quadric, $X_y=Q\cap y$ is also a quadric in
$\P2$ (of degree two) and $\O{\P2}(-X_y)=\O{\P2}(-2)$.
Plug this into the long exact sequence associated to \thetag\bandar
$$\aligned
0 & \ra H^0(\P2,\O{\P2}(m-2))\ra H^0(\P2,\O{\P2}(m))\ra
H^0(X_y,\O{X_y}(m))\\
& \ra H^1(\P2,\O{\P2}(m-2))\ra H^1(\P2,\O{\P2}(m))\ra
H^1(X_y,\O{X_y}(m)).\endaligned$$
According to Theorem B, for some suitably large $\mu_0$,
$H^q(M,\O{}(H^\mu\otimes \Cal{V}))=0$, for all $q>0,$ and $\mu>\mu_0$.
In particular, $H^1(\O{\P2}(m-2))=0$ for $m>m_0$, hence
$$h^0(\O{X_y}(m))=h^0(\O{\P2}(m))-h^0(\O{\P2}(m-2))$$
which is independent of $y$.

\head 9. The Connection and Higgs' Field \endhead

The construction of the
connection $\D{}$ in \S3 also defines a connection on $Y$.  We then have
a connection over the real slice $Y_{\Bbb R}$ of $Y$, which is a partial
compactification of $\Bbb R^{3}$.  We show that we can push $\D{}$ down
from $Y_{\R{}}$ to ${\Bbb S}^2\times\R{}$.

\subhead 9.1 The set $Y$\endsubhead
To summarise what we know about $Y$:
\roster \item Finite real points are in $Y$:
  $\{[a,b,\bar a,1]\in\C^3\subset\P3:a\in\C,\,b\in\R{}\}
=$fix$(\tau)\cap\C^3\subset Y$ since $\Cal{V}$ is trivial
on $\tau$-real sections of $T\P1$.  (In fact $t$-invariance of $\Cal{V}$
implies fix$(\sigma)\subset\{[a,b,\bar a,1]:a,b\in\C\}\subset Y$, as well.)
\item For infinite points, we know precisely that
$$\align Y\cap\Bbb{P}^{2*}_{\text{at }\infty}
&=\{\text{hyperplane sections which
contain neither }P_0\text{ nor }P_\infty\}\\
&\quad\quad(
=\frac{\{(\l_1,\l_2)\in\Cs\times\Cs\}}{(\l_1,\l_2)\sim(\l_2,\l_1)})\\
&=\{[a,b,c,0]\in\Bbb{P}^{2*}_{\text{at }\infty}: ac\ne0\}\\
&\cong \P2\setminus(\P1\vee\P1).\endalign$$\endroster
By virtue of our choice of $\C^3\subset\P3{}^*$, we know that the $\P2$ at
infinity is the set of hyperplane sections of $Q$ through the singular
point.  Furthermore, any two such hyperplane sections either have a line
in common or  meet
only at the singular point.  {}From either definition of the connection,
it's clear that the evaluation at the singular point gives a covariant
constant frame of $\Cal{W}$ over
$Y\cap\Bbb{P}^{2*}_{\text{at }\infty}$.  Since
the
standard trivialisation of $\Cal{W}$ comes from evaluation at
$P_{\lambda}$, which contains the point at $\infty$, the connection is
zero along the $\Bbb P^{2}$ at infinity.  This is
exactly the property which allows us to push the connection down to ${\Bbb
S}^2$.

\subhead 9.2 Real Points\endsubhead
We know that finite real points are in $Y$.  We must calculate the
infinite ones.
The involution $\tau$ acts on $\P{3*}$ by $\tau(a,b,c,d)=(\bar c,\bar
b,\bar a,\bar d)$, so the real points of
$Y\cap\Bbb{P}^{2*}_{\text{at }\infty}$
are $\{[a,t,\bar a,0], a\in\Bbb C^{*}, t\in\Bbb R\}\cong\Bbb R\Bbb
R^{2}\setminus \text{pt}$.  As the finite portion of $Y_{\Bbb R}$ is
$\Bbb R^{3}$, one has
$$Y_{\R{}}\cong\R{}P^3\setminus\{\pt\}.\tag9.3$$

We can similarly calculate the level set
$$\matrix Y_{\R{}}\cap\{t=0\}&=&\R3\cap\{t=0\} &\cup& (\R{}P^2\setminus
\{\pt\})\cap\{t=0\}\\
  &=&\R2 &\cup& \{[a,0,\bar a,0]:a\ne0\}\\
  &=&\R2 &\cup& {\Bbb S}^1\\
  &=& \R{}P^2\endmatrix.\tag9.4$$

\subhead 9.5 Pushing down the connection\endsubhead
Let $$\pi:\R{}P^2\to {\Bbb S}^2$$ be the map which maps the circle at
infinity in
$\R{}P^2$ to the infinite point in ${\Bbb S}^2$.

We are given a connection
on the trivial bundle over $\Bbb R P^{2}$
with a fixed trivialisation, which restricted to the
${\Bbb S}^1$ at infinity is covariant constant.  It is therefore
plausible that the connection can be pushed down.  We now check this.

Let $[x,y,w]$ be homogenous coordiates on $\R{}P^2$, with $(x/w)$ and
$(y/w)$ the usual affine coordinates on $\R2$.  Real projective
space $\R{}P^2$ has two more affine coordinate patches.  On each patch
the connection is given by connection matrices real analytic in the
affine coordinates.  In terms of the patch $(x/y),(w/y)$, the condition
that the given trivialisation is covariant constant on the circle at
infinity translates into the fact that $A_{(x/y)}=0$ when $(w/y)=0$.
Since $A$ is real analytic, we can find a second matrix such that
$A_{(x/y)}=(w/y)\widetilde{A_{(x/y)}}$.  This is analogous to the
condition
\thetag{\panda} for a connection which comes from a uniton.

Using the the relations $(x/y)=(x/w)/(y/w)$ etc., we calculate
$$\align
A_{(x/w)}&=\frac1{(y/w)^2}\widetilde{A_{(x/y)}}\\
  &=-\frac{(y/w)}{(x/w)^3}\widetilde{A_{(y/x)}}-\frac1{(x/w)^2}A_{(w/x)}\\
A_{(y/w)}&=-\frac{(x/w)}{(y/w)^3}\widetilde{A_{(x/y)}}
  -\frac1{(y/w)^2}A_{(w/y)}\\
&=\frac1{(x/w)^2}\widetilde{A_{(y/x)}}.
 \endalign$$
Discarding a neighbourhood of $((x/w),(y/w))=0$ in $\R2$, we can
integrate the squared norm of the first expression for $A_{(x/w)}$ on the
sectors
$(y-x)(y+x)>0$ and the second expression on the other sectors.  A similar
argument shows that $A_{(y/w)}$ is also finite on  $\R2$.  It follows
that the connection on $\R2$ comes from a finite-energy harmonic map or,
equivalently, that the connection pushes down to ${\Bbb S}^2$.

It will follow from the proof of $t$-invariance of $(\D{},\Phi)$ on $\R3$
that the extension of the connection from $\R2$ to ${\Bbb S}^2$ implies
the
extension of $(\D{},\Phi)$ from $\R3$ to ${\Bbb S}^2\times\R{}$.
Specifically, we
know from the discussion of the Bogomolny normalisation following
\thetag{\parana} that we can put any $t$-invariant pair $(\D{},\Phi)$
satisfying the Bogomolny equations on $\R3$ into the form $A_t=-iA_y$,
$\Phi=iA_x$, so that the finiteness of $A_x$ and $A_y$ at infinity
certainly imply the finiteness of $A_t$ and $\Phi$.

\head 10. Extra structure\endhead

\subhead 10.1 Choosing a trivialisation of $\Cal{W}\ra Y $ \endsubhead
The bundle $\Cal{V}$ is trivial over the fibre
$P_{-1}\subset\T$ with a fixed framing $\phi$.
We get a map
$\Cal{W}\ra\C^N$, defined by
$$\align\Cal{W}_y&=H^0({G}_y,\Cal{V})
  @<\text{restr.}<\cong< H^0({G}_y\cup P_{-1},\Cal{V})\\
&@>\text{restr.}>\cong> H^0(P_{-1},\Cal{V})
  @>\phi>\cong> \C^N,\tag10.2\endalign$$
where we use the canonical isomorphisms coming from restriction and
$G_{y}$ is a section of $T\Bbb P\to\Bbb P^{1}$.
This is well defined because $y$ is either a finite point (a section of
$T\P1$) and intersects $P_{-1}$ in a point, or it is infinite in which
case it meets the fibre at one point or on the whole fibre, in which case
evaluation at any point of the fibre gives the same answer, again because
$\Cal{V}$ is trivial there.
This map gives a trivialisation
$$\CD \Cal{W}@>\Psi>> \C^N\times Y\\
@VVV @VVV\\
Y @>>\text{id}> Y. \endCD$$

One value of this trivialisation is that in this framing the translation
action $\delta_t$ lifts to id$\times\delta_t$:
$$\CD
\Cal{W}_{\delta_t y} = H^0(\delta_t {G}_y,\Cal{V}) @>\eval>\cong>
\Cal{V}_{\delta_t
{G}_y\cap
P_{-1}} @<\eval<\cong< H^0(P_{-1},\Cal{V})=\C^N\\
@V\tilde \delta_t VV @V\tilde \delta_t VV @VV\tilde\delta_t=\text{id}V\\
\Cal{W}_y=H^0({G}_y,\Cal{V}) @>\eval>\cong> \Cal{V}_{{G}_y\cap P_{-1}}
@<\eval<\cong<
H^0(P_{-1},\Cal{V})=\C^N\endCD.$$
The fact that the last $\tilde\delta_t$ is the identity comes from the
fact that $\tilde\delta_t$ fixes the bundle over the section at infinity,
and hence must act trivially on sections
of the bundle over nonpolar fibres of $\T$ (over which
the bundle is trivial).
Note that $Y$ is a time-translation-independent set,
because isomorphic bundles over $\P 1$
have the same holomorphic structure.

It's not hard to see that if $\Cal{V}\to\T$ was constructed from a uniton
as in
\S\S3-7, then we have just reconstructed the original framing of
the trivial bundle
$\underline{\C}^N$.

\subhead 10.3 Time Invariance of $\D{}$ and $\Phi$\endsubhead
The connection $\D{}$ was defined by constructing covariant constant
frames along null lines in $\Bbb C^{3}$.
Consider again the covariant constant frame given by evaluation of
sections over $G_{y}$ at
$(\l_0,\eta_0)$, which defines a connection on a null plane through $y$.
Translating by $t$, we get a null plane through $y+t$, the
set of sections of
$T\P 1$ through $(\l_0,\eta_0+2\l_0 t)$.
By definition, the covariant constant frame is carried by $\tilde \delta_t$
into another covariant constant frame.

Specifically, a covariant constant frame over
$\Pi_{\l_0}$ is given as the inverse image
of a frame $f$,
$$f\in \Cal{V}(\l_0,\eta_0) @<\eval<<H^0({G}_y,\Cal{V})=\Cal{W}_y,$$
for all $y$ such that $(\lambda_{0},\eta_{0})\in G_{y}$.  Since
$$\CD
\Cal{V}(\l_0,\eta_0) @<\eval<<H^0({G}_y,\Cal{V})=\Cal{W}_y\\
@V\tilde\delta_tVV @V\tilde\delta_tVV \\
 \Cal{V}(\l_0,\eta_0-t\l_0)
@<\eval<<H^0({G}_{y+t},\Cal{V})=\Cal{W}_{y+t}\endCD$$
commutes, covariant constant frames are sent to covariant constant frames.
Hence the null connections are invariant under $\tilde\delta_t$.
Now $\D{}$ and $\Phi$ are defined in terms of these connections, so they
must perforce be invariant.
In terms of the special trivialisation, $\Psi$, the connection matrices
and the matrix representing $\Phi$ are independent of $t$.

\subhead 10.4 Reality\endsubhead
It is sufficient to show that the constructed connection and Higgs' field
satisfy our reality condition on a dense subset of $Y$.  For simplicity we
choose to work on $Y\cap\C^3$.  We can also assume that $\Cal{V}$ comes
from
$(\Cal{W},\D{},\Phi)$ which are independent of time.  It remains to show
that
they are real given a real structure on $\Cal{V}$.

We assume that the principal bundle of frames of $\Cal{V}$ comes with a
fixed antiholomorphic involution, $\tilde\sigma$, lifting $\sigma$ which
is
given in a unitary frame as $X\mapsto \overline{ X}^{t-1}$ on the fibres of
$\operatorname{fix}(\sigma)$.  As was true for $\tilde\delta_t$ and
$\Cal{W}$,
$\tilde\sigma$ induces a map on the bundle of frames of $\Cal{W}$,
$Fr(\Cal{W})\cong\Gl N$:
$$Fr(\Cal{W}_y)=H^0({G}_y,Fr(\Cal{V}))
\overset\tilde\sigma\to\leftrightarrow
H^0({G}_{\sigma(y)},Fr(\Cal{V}))=Fr(\Cal{W}_{\sigma(y)}),$$ which acts on
a section, $f(y)\in H^0({G}_y,Fr(\Cal{V}))$, by $f\mapsto\tilde\sigma\circ
f\circ\sigma$,
giving another section $f'\in H^0({G}_{\sigma(y)},Fr(\Cal{V}))$.  If
$f$ is
holomorphic, so is its image, which is the composition of one holomorphic
and two antiholomorphic maps.  The same argument as for $\tilde\delta_t$
similarly shows that $\tilde\sigma$ pulls back constant frames of the
null
connection on $\Pi$ to constant frames of the null connection on
$\sigma\Pi$, and hence the null connections back to corresponding null
connections.

In particular, if $\sigma(y)=y$, then a frame $f$ gets sent to
$\tilde\sigma(f)$ where $\tilde\sigma(f)|_{y}=\overline{(f_y)}^{t-1}$ in
any unitary basis.  So if $f$ is a covariant constant frame for $\D{}=d+A$
in
direction $X$, then $\tilde\sigma(f)$ is covariant constant for $\D{}$ in
direction $\sigma(X)$ so
$$\align A_X|_y&=\left.\dd Xf\cdot f^{-1}\right|_y\\
A_{\sigma(X)}|_y&=\left.\dd{\sigma(X)}\tilde\sigma(f)
                     \cdot\tilde\sigma(f)^{-1}\right|_y\\
&=-\left.\overline{\dd{\overline{\sigma(X)}}f\cdot f^{-1}}\right|_y\\
&=-\left.\overline{A_{\overline{\sigma(X)}}}^t\right|_y.\endalign$$
Since $\sigma(\dd x)=\dd x$, $\sigma(\dd y)=\dd y$ and $\sigma(\dd
t)=-\dd t$ (see \thetag\platypus) we have
$$A_x=-\bar A_x^t,\quad A_y=-\bar A_y^t,\quad A_t=\bar A^t_t.$$
Reality for $\Phi$  follows in a similar fashion.
To sum up, $$A_x,A_y,iA_t,i\Phi\in\u N,$$
in the chosen trivialisation of $\Cal{W}$.

\bigpagebreak
The proof of
Theorem A will be complete once we have derived the energy formula in the
next
section.

\head 11. Chern Class equals Energy \endhead
We now use the interpretation of the extended solution to
show that the second Chern class of the uniton
bundle, evaluated on the fundamental cycle of $\T$ is (a constant
multiple of) the energy of the corresponding uniton.

The Chern-Weil homomorphism tells us that the the second
Chern class of a bundle is given by integrating the trace of the
curvature squared (wedged with itself) over the base manifold for any
connection.  Define a global connection on $\T$ by gluing together
two flat connections associated to particular gauges over
$$\align U_{1}&=\{0\ne\lambda\ne\infty\}\subset\T,\text{ and}\\
U_{2}&=\{|\lambda|\ne1\}.\endalign$$
Over $U_{1}$ we take the $G_\infty$ trivialisation. Over $U_2$ we can use
the constant trivialisation of $\Bbb{C}^N\times\Bbb{R}^3\to\Bbb{R}^3$ as
follows:

The composition
$$\Bbb{R}^2\times\Bbb{P}^1\cong\{t=0\}\hookrightarrow
\Bbb{R}^3\times\Bbb{P}^1 \to T\Bbb{P}^1$$
is a (local) isomorphism onto $U_2$.  Restriction to $\{t=0\}$ of the
pullback of $\Cal{V}$ by $\Bbb{R}^3\times\Bbb{P}^1 \to T\Bbb{P}^1$ lifts
the local isomorphism to a bundle isomorphism
$$\CD
\Bbb{C}^N\times\Bbb{S}^2\times\{|\lambda|\ne1\} @>>>
\Cal{V}|_{\T\setminus\{|\lambda|=1\}} \\
@VVV  @VVV \\
\Bbb{S}^2\times\{|\lambda|\ne1\} @>>> \T\setminus\{|\lambda|=1\}.\endCD$$
We define $\nabla_2$ to be the connection which annihilates the constant
gauges of $\Bbb{C}^N\times\Bbb{S}^2\times\{|\lambda|\ne1\} \to
\Bbb{S}^2\times\{|\lambda|\ne1\}$.  (Note that this gauge cannot be
holomorphic.)

The total connection is
$$\nabla=(1-\rho)\nabla_1+\rho\nabla_2$$
where $\rho$ is a smooth function $\rho:\Bbb{P}^1\to[0,1]$ depending only
on $|\lambda|^2=\lambda\bar{\lambda}$ which is one on neighbourhoods of
$\lambda\in\{0,\infty\}$ and zero on a neighbourhood of
$\{|\lambda|=1\}$.  Since $\nabla_1$ is flat, we can integrate over $U_2$
rather than all of $\T$.  Reality allows us to integrate over one of the
components of $U_2$ and double our answer.

On $U_2$, $\nabla_2=d$ in the constant gauge over $\{t=0\}$.  Since the
$G_\infty$ trivialisation comes from the extended solution, $E_\lambda$,
lifted to the total twistor space,
$$\align \nabla_1&=d+E_\lambda^{-1}dE_\lambda\\
&\,\deq\, d+B.\endalign$$
The equations for the extended solution give
$$\align B_\lambda&=E_\lambda^{-1}\dd{\lambda}E_\lambda\\
B_{\bar{\lambda}}&=0\\
B_z&=(1+1/\lambda)A_z\\
B_{\bar{z}}&=(1+\lambda)A_{\bar{z}}.\endalign$$
So
$$\nabla=d+(1-\rho)B.$$

We calculate
$$\Omega = (\rho^{2} - \rho) B\wedge B - d\rho\wedge B,$$
and
$$\align \tr\Omega\wedge\Omega
& = (\rho^{2} - \rho)^{2}\tr B\wedge B\wedge B\wedge
B -(\rho^{2} -\rho) d\rho\wedge \tr B\wedge B\wedge B - d\rho\wedge
d\rho\wedge \tr B\wedge B \\ & =
0+3\lambda(\rho^{2}-\rho)\rho^{\prime}\tr [B_{z},B_{\bar{z}}]
B_{\lambda} dz\, d\bar{z}\, d\lambda \,d\bar{\lambda}+0 \\
 & = 3(1+\lambda)^{2}
(\rho^{2}-\rho)\rho^{\prime} \tr\dd{\bar{z}} A_{z}B_{\lambda}
dz\, d\bar{z} \, d\lambda \, d\bar{\lambda}.  \endalign$$
The first term is zero because $B$ has no $d\bar{\lambda}$ component; the
last term because
$d\rho\wedge d\rho$ is a multiple of $d\lambda\wedge d\bar{\lambda}$  and
$\tr [B_{z}, B_{\bar{z}}] = 0$.

\proclaim{Lemma 11.1} For $\lambda \in \Bbb C^{*}$ fixed,
$$\int_{z\in\Bbb C}^{}\tr\dd{\bar{z}}A_{z}B_{\lambda}dzd\bar{z} =
\lambda^{-1}\int_{z\in\Bbb C}^{}\tr A_{z}A_{\bar{z}}
dzd\bar{z}.$$ \endproclaim

\demo{Proof} Integrate by parts to obtain,
$$\int_{|z|<R}^{}\tr\dd{\bar{z}}A_{z}B_{\lambda}dzd\bar{z}
=-\int_{|z|=R}^{}\tr B_{\lambda}(A_{z}dz) - \int_{|z|<R}^{}\tr
A_{z}\dd{\bar{z}}B_{\lambda}dzd\bar{z}.$$
Since $A_{z}dz=A_{\hat{z}}d\hat{z}$ ($\hat{z}=1/z$)
and $B_{\lambda}$ is a
continuous function of $\Bbb S^{2}$, the line integral vanishes as
$R\mapsto\infty$, leaving
$$\align\int_{z\in\Bbb
C}^{}\tr\dd{\bar{z}}A_{z}B_{\lambda}dzd\bar{z}
& =-\int_{z\in\Bbb C}^{}\tr
A_{z}\dd{\bar{z}} B_{\lambda}dzd\bar{z}\\ & =-\int_{z\in\Bbb C}^{}\tr
A_{z}(\dd\lambda B_{\bar{z}} + [B_{\lambda},B_{\bar{z}}])dzd\bar{z}\ \ \ \ \ \
\ \ \ \ \ \ \ (d+B\text{ is flat}) \\ &
=-\int_{z\in\Bbb C}^{}\tr A_{z}\dd\lambda
((1+\lambda)A_{\bar{z}})dzd\bar{z} + (1+\lambda)\int_{z\in\Bbb
C}^{}\tr\, [A_{z},A_{\bar{z}}]B_{\lambda}dzd\bar{z}\\
 & =-\int_{z\in\Bbb
C}^{}\tr A_{z}A_{\bar{z}}dzd\bar{z} + (1+\lambda)\int_{z\in\Bbb
C}^{}\tr\dd{\bar{z}}A_{z}B_{\lambda}dzd\bar{z}\endalign$$
(by \thetag{1.7})
which implies the required result.\qed\enddemo

Finally, we have to integrate over $\lambda$,
$$\align\frac12\int_{}^{}\tr\Omega\wedge\Omega & = 3\int_{}^{}\tr
(1+\lambda)^{2}(\rho^{2}-\rho)\rho^{\prime}
(\lambda^{-1}A_{z}A_{\bar{z}})dzd\bar{z}d\lambda d\bar{\lambda}\\ &
=3\int_{z\in\Bbb C}^{}\tr A_{z}A_{\bar{z}}dzd\bar{z}\
 2i\int_{}^{}(\rho^{2}(r^{2})-\rho(r^{2}))\rho^{\prime}(r^{2})
 (\frac{e^{-i\theta}}{r} + 2r + re^{i\theta})drd\theta\\
& =2\pi i\int_{z\in\Bbb C}^{}\tr
A_{z}A_{\bar{z}}dzd\bar{z}.\endalign$$
So $$c_2(\Cal{V})=\frac{-1}{4\pi^2}
 \int_{\T}\tr\Omega\wedge\Omega=\frac{-i}{\pi}\int_{z\in\Bbb{C}}
 \tr A_zA_{\bar{z}}\,dz\,d\bar{z}=\frac{1}{4\pi}\operatorname{energy}(S)
 \tag11.2$$
which implies that the energy spectrum is contained in
$8\pi\Bbb{Z}$, as shown by Valli \cite{Va1}.

\head12. Ward's Construction\endhead

We are now ready to describe the link with the construction of Ward.
The value of this is that Ward's construction involves only the factoring
of a transition
matrix (i.e.\ solving the Riemann-Hilbert problem) and no differential
equations.
In addition to its metaphysical significance, this result allows us to
affirm the conjecture of Wood that unitons have rational functions in $x$
and $y$ as entries.

Let $\Cal{V}\to\T$ be a uniton bundle.
Theorem A allows us to assume that $\Cal{V}$ was constructed from a
uniton,
$S$, via $(\D{},\Phi)$, a solution to the Bogomolny equations.  We are
trying to find some intrinsic definition for $S$ on $\R3\times\P1$ which
can
be pushed down to $\T$ and interpreted as a construction for $S$.
We will in fact construct the extended
solution $E_{\lambda}$ of Uhlenbeck.

Recall the twistor fibration
$$\matrix
&&\Bbb{R}^3\times\Bbb{S}^2\\
\pi_1&\swarrow&&\searrow&\pi_2\\
\R3&&&&\T.\endmatrix$$
We defined the $G_{\infty}$ trivialisation of $\Cal{V}|_{\{0\ne
|\lambda|\ne\infty\}}$ by pulling back the extended solution
$E_{\lambda}:\Bbb S^{2}\times\Bbb C^{*}\to Gl(N)$ by
$$\CD
\Bbb R^{2}\times\Bbb R\times\Bbb C^{*}\\ @VVV\\ \Bbb S^{2}\times \Bbb
R\times\Bbb C^{*}\\ @VVV\\ \Bbb S^{2}\times\Bbb C^{*}
\endCD$$
to give a solution to
$\{\frac{\partial}{\partial t},\overline{\nabla}, \nabla_{u}-i\Phi\}$.  The
$G_{\infty}$ trivialisation is the holomorphic section of $\Cal
V |_{\{0\ne\lambda\ne\infty\}}$ represented by this solution.  On the
other hand, a frame of $\Cal W$ over a point, $y\in\Bbb R^{3}$, lifts to
a trivialisation of $\pi_{1}^{*}\Cal W$ restricted to the complex line
$y\times\Bbb{S}^2$.
Since this trivialisation is, by definition, in the kernel of
$\overline{\nabla}$, it pushes down to a trivialisation of $\Cal V
|_{G_{y}}$.

Now we use basic holomorphic geometry to reconstruct the extended
solution.  Since holomorphic functions on $\Bbb P^{1}$ are constant,
holomorphic frames of $\Bbb C^{N}\times\Bbb P^{1}\to\Bbb P^{1}$ are
determined by evaluation at a point.  Since $\Cal V |_{G_{\infty}}$ and
$\Cal V |_{P_{\lambda}}$, $\lambda\ne 0, \infty$, are trivial, the
$G_{\infty}$ trivialisation is intrinsically defined (by fixing the
framing $\phi$), as is the trivialisation of $\Cal V |_{G_{y}}$ which
agrees with the framing $\phi$ at $G_{y}\cap P_{-1}$.  The change of
gauge $\Cal{E}(\lambda, G_{y})$ which relates these two trivialisations
at $(\lambda,\eta)$, expressing the $G_{\infty}$
trivialisation in terms of the $G_{y}$ trivialisation of $\Cal V
|_{G_{y}\cap P_{\lambda}}$ is thus well defined.
(Note that at $\lambda=-1$, $\Cal{E}(-1,G_y)=\one$.)
Lifting to the total twistor space, $\Cal{E}(\lambda ,G_y)$, is
the gauge transformation from the $\pi_{2}$-pulled-back $G_{\infty}$
trivialisation to the $\pi_1$-pulled-back standard gauge, and so
$\Cal{E}(\lambda ,G_y)$ is the extended solution, $E_{\lambda}$.

Although we have used the total twistor fibration to interpret the
extended solution, we have only {\it required} the framing and the trivial
holomorphic structure of $\Cal V$ restricted to certain lines.
We can thus express $E_{\lambda}$ intrinsically (without lifting to
the total space) as
the `monodromy' around the cycle of $\P1$'s in \figh.
\hpict  

To make this precise, what we are calling a `monodromy' is actually the
failure to commute of a cycle of homomorphisms given by the restriction
map:
$$\CD
\Cal{V}_{\l,\infty} @<\text{restr}<< H^0({G}_\infty,\Cal{V})
@>\text{restr}>>
  \Cal{V}_{-1,\infty}\\
@A\text{restr}AA @. @AA\text{restr}A\\
H^0(P_\l,\Cal{V}) @. @. H^0(P_{-1},\Cal{V})\\
@V\text{restr}VV @. @VV\text{restr}V\\
\Cal{V}_{(\l,z/2-\l t-\l^2\z/2)} @<\text{restr}<<
H^0({G}_{(z,\z,t)},\Cal{V})
  @>\text{restr}>> \Cal{V}_{(-1,z/2+ t-\z/2)}.\endCD\tag12.1$$
The `monodromy' is independent of the choice of initial value, up to
conjugation, as one would expect, since a change of framing of the bundle
acts by conjugation on the uniton.  We fix it by computing the
`monodromy' of the fixed frame $\phi\in H^0(P_{-1},Fr(\Cal{V}))$.

\subhead12.2 Transition Functions \endsubhead
Ward's construction assumes the bundle is given by a transition matrix, so
consider the covering of $\T$ given by
$$\aligned
U&=\{\l\in\C,\eta\in\C\},\\
\hat U&=\{\lh\in\C,\eh\in\C\},\\
U'&=\{\l\in\C,\ep\in\C\},\\
\hat U'&=\{\lh\in\C,\ehp\in\C\}.\endaligned\tag12.3$$
The bundle $\Cal{V}$ is determined by transition matrices $T,\hat T,T'$
which
map fixed frames of $\Cal{V}$ over $U$ to $\hat U$, over $\hat U$ to $\hat
U'$,
and over $U$ to  $U'$ respectively.
Because $\Cal{V}$ has certain triviality properties, we can choose the
fixed
frames such that
$$\align
\hat TT{T'}^{-1}|_{{G}_\infty}&=\one, \\
T'|_{P_1}&=\one, \text{ and}\tag\ibix\\
\hat T|_{P_{-1}}&=\one.\endalign$$

If the bundle is trivial when restricted to a complex line ($\P1$), then a
framing above a point of the line extends uniquely to a nonvanishing frame
on
the line, because in this case evaluation
$$H^0(\P1,\underline{\C}^N)@>\text{eval}>>\underline{\C}^N_{\text{point}}
=\C^N$$ is an isomorphism.
We will think of this as defining a parallel translation within the line.

In terms of these frames parallel translation from a point on $P_1$ (in
terms of the $U$ frame) to a point on $P_{-1}$ (in terms of the $\hat U$
frame) along $P_1\cup {G}_\infty\cup P_{-1}$ is given by $\one$.
Since the bundle is trivial above real sections, we can split
$T$, i.e.\ find analytic functions $H:\{(z,\z,t)\in\R3,\l\in\C\}\to\Gl N$ and
$\hat H:\{(z,\z,t)\in\R3,\l\in(\P1\setminus\{0\})\}\to\Gl N$, such that
$$TH_\l(z,\z,t)=\hat H_\l(z,\z,t).$$
Parallel translation from $P_1\cap {G}_y$ (in terms of the $\hat U$ frame) to
$P_{-1}\cap {G}_y$ ($U$ frame) along ${G}_y$ is given by $$
H_{1}(y)\hat H_{-1}(y)^{-1}\ \(=E_1(y)=S(y)\),\tag\emu$$
which gives the same formula for the uniton as in \cite{\Ward,\,18}
(up to a change of sign caused by a difference in basing conventions).
One must verify that this doesn't depend on the choice of splitting.

Finally, we note that Ward actually takes {\it two\/} framings, along
$P_1$ and $P_{-1}$.  One framing would be equivalent to the restrictions
\thetag\ibix.  By taking two framings he does away with the basing
condition.  This is an important point if one wants to choose a different
type of basing condition (other than $S(\infty)=\one)$, to encode
Grassmannian solutions for example.

\head 13. Example of a $U(3)$ uniton.\endhead

Over a ruled surface, holomorphic bundles framed along a section of the ruling
which are trivial on the generic
line of the ruling are determined by their \lq jumps\rq\
 (i.e.\ the structure of the bundle on a neighbourhood of the
 \lq jumping\rq\ lines on which the bundle is not trivial).
(See \cite{Hu2}, \cite{HuMi}.)
If $\eta$ is a coordinate on the lines and $\lambda$ parametrises the
lines of the ruling, with $\lambda=0$ a jumping line,
then the \lq jump\rq\ is determined by a transition matrix
$J(\lambda,\eta,\eta^{-1})$ from $\{\eta\ne\infty\}$ to $\{\eta\ne0\}$.
It can then be \lq glued\rq\ to the bundle over $\lambda\ne0$
(which is trivialised by the framing of the bundle over a section of
the ruling) by an additional matrix $M(\lambda)$.
In our case, the $G_{\infty}$ trivialisation
is the \lq generic\rq\ trivialisation away from $\lambda=0$.

Assuming a \lq unitary\rq\ choice of frame, $\phi$,
the $G_{\infty}$ trivialisation is fixed by
the real structure
($\tilde\sigma$) and the real structure exchanges the  jumps
 at $P_{0}$ and
$P_{\infty}$. So the uniton bundle is determined by
the jump at $P_0$, $J(\lambda ,\eta,\eta^{-1})$, and its
 \lq gluing\rq, $M(\lambda)$, alone.
As an example we will show how to construct a
$U(3)$ uniton from such
 a jump.

Let
$$\pmatrix \eta^{2}&2\lambda\eta&\lambda^{2}\\
   0&1&{\frac {\lambda}{\eta}}\\
   0&0&\eta^{-2} \endpmatrix\tag 13.1$$
be the transition matrix for the jumping line $P_0$, and let this jump
be glued in by $M(\lambda)=\one$ over the
section at infinity to the trivial
bundle on the union of the nonpolar fibres.  For the
bundle to admit a lift of the real structure, the transition matrix
at $P_\infty$ must
be conjugate to \thetag{13.1}, i.e.\
$$\pmatrix \hat{\eta}^{-2}&0&0\\
2\frac {\hat{\lambda}}{\hat{\eta}}&1&0\\
\hat{\lambda}^{2}& \hat{\lambda}\hat{\eta} &\hat{\eta}^{2}\endpmatrix.$$
We have thus defined a bundle over $\widetilde{T\Bbb P}{}^1$, with the
required reality and triviality over fibres and the section at infinity.
We need to check that it also has time translation and that it
is trivial when restricted to real sections of $T\Bbb P^1$.

\subhead 13.2 ${G}_\infty$ trivialisation\endsubhead
Near $\lambda=0$, the $G_\infty$ trivialisation can be written as
$$\pmatrix  0&0&\lambda^{2}\\0&-1&-\lambda\eta
\\\lambda^{-2}&{2\frac {\eta}{\lambda}}&\eta^{2}\endpmatrix
\qquad \text{and} \qquad
\pmatrix 1&0&0 \\
\frac {1}{\lambda\eta}&1&0 \\
\frac {1}{\lambda^{2}\eta^{2}}&\frac {2}{\lambda\eta}&1 \endpmatrix $$
in the
$\{\eta\ne\infty\}$
and $\{\eta\ne0\}$ trivialisations, respectively.
(These are the  sections' trivialisations
of \cite{Hu2}.)

Similarly, near $\hat{\lambda}=0$ the $G_\infty$ trivialisation
can be written as
$$\pmatrix
 \hat{\eta}^{2}&{\frac {\hat{\eta}}{\hat{\lambda}}}&
\hat{\lambda}^{-2}\\-2 \hat{\lambda}\hat{\eta}&-1&0
\\\hat{\lambda}^{2}&0&0 \endpmatrix
\qquad \text{and} \qquad
\pmatrix
 1&{\frac {1}{\hat{\lambda}\hat{\eta}}}&{\frac {1
}{\hat{\eta}^{2}\hat{\lambda}^{2}}}\\0&1&{\frac {2}{\hat{\lambda}
\hat{\eta}}}\\0&0&1 \endpmatrix.$$
We can use these trivialisations to check that our bundle admits time
translation and is trivial on real sections.

\subhead 13.3 Time Translation \endsubhead
Time translation is determined by
the requirement that it fixes the bundle above the infinity section.
This means that it sends the framing
$$\pmatrix  0&0&\lambda^{2}\\ 0&-1&-\lambda\eta
\\\lambda^{-2}&{2\frac {\eta}{\lambda}}&\eta^{2}\endpmatrix
 \qquad \text{to} \qquad
\pmatrix  0&0&\lambda^{2}\\ 0&-1&-\lambda(\eta+\lambda t)
\\\lambda^{-2}&{2\frac {(\eta+\lambda t)}{\lambda}}&(\eta+\lambda t)^{2}
\endpmatrix$$
so it is
$$\pmatrix
 1&0&0\\-t&1&0
\\t^{2}&-2t&1\endpmatrix$$
in terms of the first frame.  Since this is continuous at $\lambda=0$, time
translation extends as a continuous map to the jumping line.

Notice that on the big open set (and at $\eta=\infty$) the
time translation acts trivially (in the $G_\infty$ trivialisation).

\subhead 13.4 Real Triviality \endsubhead  Over the big open set
${0 \ne \lambda \ne \infty}$, different choices of global frames are
related by $\operatorname{Gl}(3)$-valued
functions of
$(\lambda ,\eta)$.  We want functions which,
when restricted to real sections,
extend to $\lambda \in \Bbb P^{1}$.  To test the extension we need to
change to the appropriate frame over the open sets covering $P_{0}$ and
$P_{\infty}$.  The transition matrix restricted to ${\eta=\infty}$ is the
gluing map $(\Bbb I)$, but {\it close to} $\eta=\infty$ it is given more
generally by the sections\rq\ trivialisations over $\{\eta\ne0\}$:
$$\pmatrix 1&0&0\\
\frac {1}{\lambda\eta}&1&0\\
\frac {1}{\lambda^{2}\eta^{2}}&\frac {2}{\lambda\eta}&1 \endpmatrix
_{P_{0}}
\qquad \text{and} \qquad \pmatrix
 1&{\frac {1}{\hat{\lambda}\hat{\eta}}}&{\frac {1
}{\hat{\eta}^{2}\hat{\lambda}^{2}}}\\
0&1&{\frac {2}{\hat{\lambda}\hat{\eta}}}\\0&0&1
\endpmatrix _{P_{\infty}}$$
respectively.  Near $\eta=0$, we have to use the transition
matrices to switch frames again, but this is the same as multiplying the
original frame by the first sections\rq\ trivialisation there.

After some experimenting, we find that the frame
$$\pmatrix
z^{2}&-\left (2z+z^{2}\bar{z}\right )\lambda&\lambda
^{2}\\-{\frac {z}{\lambda}}&1-z^{2}\bar{z}^{2}&2\bar{z}\lambda\\
\lambda^{-2}&{\frac {2z\bar{z}^{2}+\bar{z}}{\lambda}}&
\bar{z}^{2} \endpmatrix\tag13.5$$
extends.

For example
$$\align&\pmatrix 0&0&\lambda^{2}
\\0&-1&-\lambda\eta\\\lambda^{-2}&{\frac {2\eta}{\lambda
}}&\eta^{2}\endpmatrix
\pmatrix  z^{2}&-\left (
2z+z^{2}\bar{z}\right )\lambda&\lambda^{2}\\-{\frac {z}{\lambda}}&1
-z^{2}\bar{z}^{2}&2\bar{z}\lambda\\\lambda^{-2}&{\frac {2
z\bar{z}^{2}+\bar{z}}{\lambda}}&\bar{z}^{2}\endpmatrix\\
&=\pmatrix 1&\lambda\bar{z}\left (2\bar{z}z+1\right )&\lambda^{2
}\bar{z}^{2}\\-{\frac {-z+\eta}{\lambda}}&z^{2}\bar{z}^{2}-
1-2\eta z\bar{z}^{2}-\eta\bar{z}&-2\bar{z}\lambda-\lambda\eta\bar{z}^{2}
\\{\frac {z^{2}-2\eta z+\eta^{2}}{\lambda^{2}}}&{\frac {-2z-z
^{2}\bar{z}-2\eta z^{2}\bar{z}^{2}+2\eta+2\eta^{2}\bar{z}^{2}z+\eta^{2}
\bar{z}}{\lambda}}&1+4\eta\bar{z}+\eta^{2}\bar{z}^{2}\endpmatrix,\endalign
$$
which is seen to be holomorphic in $\lambda$ on the real section
$\eta=z-\lambda^2\bar{z}$.

\subhead 13.6 Formula \endsubhead
The nice thing is that the frames \thetag{13.5} (in terms of the
generic gauge) give us the transition from the
${G}_{(z,\overline{z},0)}$
gauge into the $P_{\lambda}/{G}_{\infty} (\lambda \in \Bbb C^{*})$
gauges,
and we can thus calculate the uniton by  Ward's method to be
$$\align E_\lambda&=(\text{framing at } \lambda)(\text{framing at }
-1)^{-1}\\
&=\frac{1}{\left (1+z^{2}\bar{z}^{2}+
\bar{z} z\right )\left(4 \bar{z} z+ z^{2}\bar{z}^{2}+1\right)}\endalign$$
$$
\{\lambda^{-2}\pmatrix
0&0 &0 \\
0&0&0\\
4 \bar{z}^3 z+ z^{2}\bar{z}^{4}+ \bar{z}^{2} & 4 \bar{z}^2
z+z^{2}\bar{z}^{3}
+\bar{z} &  z^{2}\bar{z}^{2}+4 \bar{z} z+ 1\\ \endpmatrix\right.
$$ $$
+\lambda^{-1}\pmatrix
0& 0& 0\\
-z^{3}\bar{z}^{4}-4 z^{2}\bar{z}^{3}-\bar{z}^2 z&
-z^{3}\bar{z}^{3}-4 z^{2}\bar{z}^{2}-\bar{z} z &
-z-4\bar{z} z^2-z^{3}\bar{z}^{2} \\
 2 z^{2}\bar{z}^{4}+5 z\bar{z}^3+ 2 \bar{z}^{2}&
- 2 z^{3}\bar{z}^{4}-z^{2}\bar{z}^{3}+2 z\bar{z}^2+\bar{z} &
 -4 z^{3}\bar{z}^{3}-4 z^{2}\bar{z}^{2}-z\bar{z} \\
\endpmatrix
$$ $$
+\lambda^{0}\pmatrix
z^{4}\bar{z}^{4}+4 z^{3}\bar{z}^{3}+z^{2}\bar{z}^{2}&
z^{4}\bar{z}^{3}+\bar{z} z^2 + 4 z^{3}\bar{z}^{2}&
4\bar{z} z^3+ z^{4}\bar{z}^{2}+ z^2\\
- z^{3}\bar{z}^{4}-2 z^{2}\bar{z}^{3}+z\bar{z}^2+2\bar{z}&
z^{4}\bar{z}^{4}-2 z^{2}\bar{z}^{2}+1  &
+ 2 z^{4}\bar{z}^{3}+z^{3}\bar{z}^{2}-2 z^2\bar{z}-z \\
\bar{z}^{2}+ z^{2}\bar{z}^{4}+\bar{z}^3 z&
- 2 \bar{z}^2 z- 2 \bar{z}^{4}z^{3}- 2 z^{2}\bar{z}^{3} &
z^{4}\bar{z}^{4}+\bar{z}^{3}z^{3}+z^{2}\bar{z}^{2} \\
\endpmatrix
$$ $$
+
\lambda^{1}\pmatrix
-z^{3}\bar{z}^{3}-4 z\bar{z}-4 z^{2}\bar{z}^{2}&
-2 z-z^2\bar{z} +2 z^{3}\bar{z}^{2}+z^{4}\bar{z}^{3}&
2 z^{4}\bar{z}^{2}+ 5 z^3\bar{z}+ 2  z^2 \\
+2 z^{2}\bar{z}^{3}+2 \bar{z}+2 \bar{z}^2 z&
-4 \bar{z}^{3}z^{3}-4 z^{2}\bar{z}^{2}-4 \bar{z} z &
+ 2 \bar{z} z^2+2 \bar{z}^{3}z^{4}+ 2 z^{3}\bar{z}^{2} \\
0& 0& 0 \\
\endpmatrix
$$ $$\left.+
\lambda^{2}\pmatrix
\bar{z} z+ 1+ z^{2}\bar{z}^{2}&
-2 z^{3}\bar{z}^{2}-2z-2\bar{z} z^2 &
z^2 + \bar{z} z^3 + z^{4}\bar{z}^{2}\\
0&0 & 0\\
0&0&0 \\
\endpmatrix \}
$$

As pointed out by Francis Burstall, it is easy to test for holomorphic
(one uniton) solutions.  Assume $S$ were a one uniton, corresponding to a
projection $\pi$, then its extended solution would have the form
$$E_\lambda=g(\lambda)(\pi-\lambda\pi^\perp)$$
and
$$\align { (E_\lambda |_{z=0})^{-1}E_\lambda}
&=(\pi-\lambda\pi^\perp)^{-1}|_{z=0}
  (\pi-\lambda\pi^\perp)\\
&=(\pi|_{z=0}-\lambda^{-1}\pi^\perp|_{z=0})(\pi-\lambda\pi^\perp)
\endalign$$
would contain no powers $\lambda^{-2}$ or $\lambda^2$.  Since $E_\lambda$
as constructed above
does give such terms, our example must be a  two uniton.

\head 14. Wood's Conjecture\endhead

Noting that all known examples of unitons were matrices of functions
rational in $x$ and $y$ (equivalently $z$ and $\z$),
Wood conjectured that this is always the case (\cite{\Wooda}).
While it is true that the Bogomolny solution
$(\D,\Phi)$ constructed from a uniton bundle is algebraic
we do not know that the integration
$$S^{-1}dS=2\(A_zdz+A_\z d\z\)$$
preserves rationality.
Continuing $S$ analytically, or equivalently, integrating $A$, we
can't even rule out multivaluedness if $A$ is holomorphic on
nonsimply-connected domains.

The concrete expression \thetag\emu, however, shows that $S$ extends
to $Y\cap\C^3$
(a Zariski open set), and using \thetag\emu\ and the jumping-line normal
form for transition matrices (see \cite{Hu2} and \cite{\new} for
proofs) we can prove Corollary B.

\remark{Remark 14.1}  In \cite\thesis\ we give an obviously rational,
explicit formula for $S$ in terms of monad data of a special type.  If
this could be extended to general monad data, and hence general unitons,
this proof would not be necessary. \endremark

\demo{Proof of Corollary B} Since the solution is $t$-invariant, we can
ignore the third dimension.

We want to show $E_\l(z,\z)=\hat H_{\lambda}(z,\z)H_{-1}(z,\z)^{-1}$
 is a rational $\gl
N$-valued function on $\{(z,\z)\in\P1\times\P1\}$.  A {\it function\/} is
rational \Iff it is meromorphic \Iff it is meromorphic when restricted to
the sets of a covering of $\P1\times\P1$, and a function is meromorphic
\Iff its only singularities are poles.  Thus we can answer a global
question with a local answer.  This is an example of Serre's GAGA
principle (\cite{Se}),
which says that analytic objects on compact varieties
are algebraic.

Consider the family of open sets $$\{U_{z_0}=\{(z,w)\in\P1\times\P1:z\ne
z_0,w\ne\z_0\}_{z_0\in\P1}\}.$$  Any three sets cover $\P1\times\P1$.
Since the uniton equations are conformally invariant,
the space of conformal maps $\operatorname{conf}(\Bbb{S}^2,\Bbb{S}^2)$
acts on the space of unitons.  Our family of open sets is
the orbit of $\Bbb{C}\subset\Bbb{P}^1\cong\Bbb{S}^2$
under conformal changes $z\mapsto1/(z-a)$.  Each such change of
coordinates produces a different uniton and a different uniton
bundle.  Rather than showing that a single uniton is meromorphic
on enough sets of the family, it is easier to show that all its
transformations under $\operatorname{conf}(\Bbb{S}^2,\Bbb{S}^2)$
 are meromorphic on $U_\infty$, the basic open set which
corresponds to our choice of coordinates.  It is sufficient to
do this for an arbitrary uniton.  Let $S$ be such a uniton.

The expression \thetag\emu\ defines $S$ on $\R2$, but extends just as well
to $\C^2$ with potential singularities at the jumping lines.
To see that they are poles, pull back the transition matrix $T$ by
the sections' map
$$\aligned \varphi:\C^2\times\P1&\to T\P1\\
(z,\z)\times(\l)&\to(\l,\eta=z-\z\l^2).\endaligned$$
If $(z,\z,0)$ represents a jumping line of type $(k_1\le k_2\le\cdots\le
k_N)$
(i.e.\ $\Cal{V}|_{{G}_{(z,\z,t)}}\cong\O{}(k_1)\oplus\cdots\O{}(k_N)$),
then we can make a holomorphic change of frame on some neighbourhood of
the point so that $T$ has the form
$$T=\p \lambda^{-k_1}\\
&\ddots& &p^i_j\\
& &\ddots\\
&0&&\ddots\\
&&&&\lambda^{-k_N}\pp, \quad p^i_j
=\sum_{a=-k_j+1}^{-k_i-1}p^i_{ja}(z,\z)\l^a$$
where $p^i_{ja}$ are holomorphic functions.
A section of $\Cal{V}|_y$ is given by $(u^1,\dots,u^N)^t$,
$u^j=\sum_{a=0}^\infty u^j_a\l^a$ such that $$T(y,\l)\p u^1\\\vdots\\
u^N\pp\text{ is holomorphic in }1/\l,$$
which puts conditions on $\{u^j_a\}$.
Expanding the columns of $T\cdot u$ in Laurent series in $\l$, the
conditions come from the coefficients of positive terms in $\l$ ($\l^i$,
$i>0$), which we can see are linear in $u^j_a$.

The result is a matrix equation for $U=(u^j_a)$ of the form
$$M(p^i_{ja})U=0$$
where the coefficients of $M$ are linear
polynomials in the $p^{i}_{ja}$.  Since the $p^{i}_{ja}$ are holomorphic
in $z$ and $\bar{z}$, $M$ is a holomorphic function of $y$.
The  matrix $M$ can be reduced by
eliminating coefficients $u^{j}_{a}$ which are completely determined by
other coefficients.  The reduced matrix $\Gamma(y)$ whose elements are
now polynomials in the $p^{i}_{ja}$ has the property that
$$h^{0}(G_{y},\Cal{V})=\operatorname{corank}\Gamma(y),$$
and $\Gamma(y)$ has maximal rank when $\Cal{V}|_{G_{y}}$ is trivial.  We
can parametrise the kernel of $\Gamma(y)$ in a punctured neighbourhood of
a jump by rational polynomials in $p^{i}_{ja}$.  ($\Gamma(y)$ is
polynomial, but to parametrise the kernel we must invert a maximal rank
submatrix of $\Gamma(y)$, introducing poles at the jump where the rank of
$\Gamma$ drops.)  {}From this parametrisation of the kernel we
can reconstruct the holomorphic frames of $\Cal{V}|_{G_{y}}$ (as a
function of $y$) by rational
polynomials in the $p^{i}_{ja}$, i.e.\ with possible poles at the jumping
line.  Since the extended solution is constructed from these frames by
algebraic operations, it is also meromorphic.\qed
\enddemo

\remark{Acknowledgement}  For advice on matters of mathematics and
mathematical exposition, I am indebted to Professor J.\ C.\ Hurtubise.
\endremark

\newpage

\head Remark on Chern-Simons
  classes not appearing in other versions of paper\endhead
As pointed out to the author by K.\  Guruprasad,
it is known that the energy of harmonic maps of surfaces into
groups can be calculated as the evaluation of a Chern-Simons
form.  The following example of this may be of interest.

In section 11 we
remarked that as a limiting case, the curvature of the combined
connection can be localised at the jumping lines.  This reflects the fact
that the Chern class of a bundle over a ruled surface can be calculated
locally at the \lq\lq jumping lines\rq\rq (\cite{Hu2}).
 At the other extreme, as we reduce the support of $(1-\rho)$ by taking
a family of  partition functions depending on
$\epsilon$ such that
$$\rho(t)=\{\matrix 1 & t\in[0,1-2\epsilon]\\
\text{strictly decreasing} & t\in[1-2\epsilon,1-\epsilon]\\
0 &t\in[1-\epsilon,1]\endmatrix\right.$$
and considering the limit as $\epsilon\mapsto 0$,
the Chern-Weil integral becomes the
evaluation of the Chern-Simons form of the connection $\Cal D_{\lambda}$
on the trivial bundle over $\Bbb S^{2}\times\Bbb S^{1}$.
(Strictly speaking, $\Cal D_{\lambda}$ is a family of connections on
$\Bbb S^{2}$ but we may extend it uniquely, up to  changes
of gauge holomorphic in $\lambda$.  We can fix the extention
 by fixing an extended  solution $E_{\lambda}$.)

Use coordinates $\lambda =r\exp(i\theta)$.  Then $(z,\bar{z},\theta)$ are
coordinates on $\Bbb S^{2}\times\Bbb S^{1}$, and we have a smooth family
of connections $\widetilde{B}(s)$ given by restricting $(\Cal W,B)$
to $\{r=s>0\}\cong
\Bbb S^{2}\times\Bbb S^{1}$.  Since $\rho$ is monotonic, and the
integrand is zero when $\rho$ fails to be bijective, we can change
variables to $t=\rho(r)$ and obtain
$$\align-\int_{r=0}^{1}&(\rho^{2}-\rho)\rho^{\prime}r\,dr \wedge
\text{tr}B\wedge B\wedge
B \\
& =-\int_{t=1}^{0}(t^{2}-t)\text{tr}\widetilde B(\rho^{-1}(t))\wedge
\widetilde B(\rho^{-1}(t))\wedge \widetilde B(\rho^{-1}(t))dt\\ &
\to\int_{t=0}^{1}(t^{2}-t)\text{tr}\widetilde B(1)\wedge \widetilde B(1)
\wedge \widetilde B(1)dt\endalign$$
as $\epsilon\mapsto 0$ which is the Chern-Simons class of the
{\it flat} connection $(\Cal W|_{\{|\lambda|=1\}},\Cal D_{\lambda})$
over $\Bbb S^2\times\Bbb S^1$ (\cite{ChSi}).

\enddocument